\input amstex
\documentstyle{amsppt}
\input epsf
\input texdraw

\magnification1200
\hsize13cm
\vsize19cm

\TagsOnRight

\define\twoline#1#2{\line{\hfill{\smc #1}\hfill{\smc #2}\hfill}}
\define\twolinetwo#1#2{\line{{\smc #1}\hfill{\smc #2}}}
\define\twolinethree#1#2{\line{\phantom{poco}{\smc #1}\hfill{\smc
#2}\phantom{poco}}}

\def\mypic#1{\epsffile{#1}}




\def\fba{2.1}
\def\fbb{2.2}

\def\FB{3.1}
\def\FC{3.2}
\def\FD{3.3}
\def\FE{3.4}

\def\fea{5.1}


\def\eaa{1.1}

\def\eba{2.1}
\def\ebb{2.2}
\def\ebc{2.3}
\def\AA{2.4}
\def\AB{2.5}
\def\AC{2.6}

\def\AD{3.1}
\def\AE{3.2}
\def\AF{3.3}
\def\AG{3.4}
\def\AH{3.5}
\def\AHa{3.6}
\def\AI{3.7}
\def\AJ{3.8}
\def\AK{3.9}
\def\AL{4.1}
\def\AM{4.2}

\def\eea{5.1}
\def\eeb{5.2}
\def\eec{5.3}


\def\tba{2.1}

\def\tea{5.1}

\def\TA{\tba}

\def\TC{3.1}
\def\TD{3.2}
\def\TE{3.3}
\def\TF{3.4}
\def\TG{4.1}

\def\AndrAK{1}
\def\AndrAW{2}
\def\BresAO{3}
\def\CiucAB{4}
\def\CiucNnineNten{5}
\def\CiKrAA{6}
\def\CiKrfiveandahalf{7}
\def\CiKrAF{8}
\def\DT{9}
\def\GeViAB{10}
\def\GordAB{11}
\def\KK{12}
\def\KaKZAC{13}
\def\Kuo{14}
\def\KupNnine{15}
\def\KupExploration{16}
\def\LindAA{17}
\def\MacMAA{18}
\def\OkadAA{19}
\def\Proc{20}
\def\StanSymmClass{21}
\def\StemAE{22}
\def\StemAG{23}

\def\({\left(}
\def\){\right)}

\def\sgn{\operatorname{sgn}}
\def\Pf{\operatorname{Pf}}
\def\M{\operatorname{M}}
\def\P{{\Cal P}}
\def\GF{\operatorname{GF}}

\def\Z{\Bbb Z}


\def\RhombusA{\bsegment
  \rlvec(0.866025403784439 .5) \rlvec(0.866025403784439 -.5) 
  \rlvec(-0.866025403784439 -.5) \rlvec(-0.866025403784439 .5) 
  \savepos(0.866025403784439 -.5)(*ex *ey)
        \esegment
  \move(*ex *ey)
        }
\def\RhombusB{\bsegment
  \rlvec(0.866025403784439 .5) \rlvec(0 -1)
  \rlvec(-0.866025403784439 -.5) \rlvec(0 1) 
  \savepos(0 -1)(*ex *ey)
        \esegment
  \move(*ex *ey)
        }
\def\RhombusC{\bsegment
  \rlvec(0.866025403784439 -.5) \rlvec(0 -1)
  \rlvec(-0.866025403784439 .5) \rlvec(0 1) 
  \savepos(0.866025403784439 -.5)(*ex *ey)
        \esegment
  \move(*ex *ey)
        }

\def\hdSchritt{\bsegment
  \lpatt(.05 .13)
  \rlvec(0.866025403784439 -.5) 
  \savepos(0.866025403784439 -.5)(*ex *ey)
        \esegment
  \move(*ex *ey)
        }

\def\odaSchritt{\bsegment
  \lpatt(.05 .13)
  \rlvec(0.866025403784439 .5) 
  \savepos(0.866025403784439 .5)(*ex *ey)
        \esegment
  \move(*ex *ey)
        }

\def\ringerl(#1 #2){\move(#1 #2)\fcir f:0 r:.15}
\def\Ringerl(#1 #2){\move(#1 #2)\lcir r:.19}

\catcode`\@=11
\font\tenln    = line10
\font\tenlnw   = linew10
\font\tencirc  = lcircle10    
\font\tencircw = lcirclew10   

\newskip\Einheit \Einheit=0.5cm
\newcount\xcoord \newcount\ycoord
\newdimen\xdim \newdimen\ydim \newdimen\PfadD@cke \newdimen\Pfadd@cke

\newcount\@tempcnta
\newcount\@tempcntb

\newdimen\@tempdima
\newdimen\@tempdimb

\newdimen\@wholewidth
\newdimen\@halfwidth
\newdimen\unitlength \unitlength =1pt

\newcount\@xarg
\newcount\@yarg
\newcount\@yyarg
\newbox\@linechar
\newbox\@tempboxa
\newdimen\@linelen
\newdimen\@clnwd
\newdimen\@clnht

\newif\if@negarg

\def\@whilenoop#1{}
\def\@whiledim#1\do #2{\ifdim #1\relax#2\@iwhiledim{#1\relax#2}\fi}
\def\@iwhiledim#1{\ifdim #1\let\@nextwhile=\@iwhiledim
        \else\let\@nextwhile=\@whilenoop\fi\@nextwhile{#1}}

\def\@whileswnoop#1\fi{}
\def\@whilesw#1\fi#2{#1#2\@iwhilesw{#1#2}\fi\fi}
\def\@iwhilesw#1\fi{#1\let\@nextwhile=\@iwhilesw
         \else\let\@nextwhile=\@whileswnoop\fi\@nextwhile{#1}\fi}

\def\thinlines{\let\@linefnt\tenln \let\@circlefnt\tencirc
  \@wholewidth\fontdimen8\tenln \@halfwidth .5\@wholewidth}
\def\thicklines{\let\@linefnt\tenlnw \let\@circlefnt\tencircw
  \@wholewidth\fontdimen8\tenlnw \@halfwidth .5\@wholewidth}
\thinlines

\PfadD@cke1pt \Pfadd@cke0.5pt
\def\PfadDicke#1{\PfadD@cke#1 \divide\PfadD@cke by2 \Pfadd@cke\PfadD@cke \multiply\PfadD@cke by2}
\long\def\LOOP#1\REPEAT{\def\BODY{#1}\ITERATE}
\def\ITERATE{\BODY \let\next\ITERATE \else\let\next\relax\fi \next}
\let\REPEAT=\fi
\def\Punkt{\hbox{\raise-2pt\hbox to0pt{\hss$\ssize\bullet$\hss}}}
\def\DuennPunkt(#1,#2){\unskip
  \raise#2 \Einheit\hbox to0pt{\hskip#1 \Einheit
          \raise-2.5pt\hbox to0pt{\hss$\bullet$\hss}\hss}}
\def\NormalPunkt(#1,#2){\unskip
  \raise#2 \Einheit\hbox to0pt{\hskip#1 \Einheit
          \raise-3pt\hbox to0pt{\hss\twelvepoint$\bullet$\hss}\hss}}
\def\DickPunkt(#1,#2){\unskip
  \raise#2 \Einheit\hbox to0pt{\hskip#1 \Einheit
          \raise-4pt\hbox to0pt{\hss\fourteenpoint$\bullet$\hss}\hss}}
\def\Kreis(#1,#2){\unskip
  \raise#2 \Einheit\hbox to0pt{\hskip#1 \Einheit
          \raise-4pt\hbox to0pt{\hss\fourteenpoint$\circ$\hss}\hss}}

\catcode`\@=11
\newif\if@ovt
\newif\if@ovb
\newif\if@ovl
\newif\if@ovr
\newdimen\@ovxx
\newdimen\@ovyy
\newdimen\@ovdx
\newdimen\@ovdy
\newdimen\@ovro
\newdimen\@ovri
\newcount\@tempcnta
\newcount\@tempcntb
\newdimen\@tempdima
\newdimen\@tempdimb
\newdimen\@wholewidth
\newdimen\@halfwidth
\newdimen\unitlength \unitlength =1pt

\def\thinlines{\let\@linefnt\tenln \let\@circlefnt\tencirc
  \@wholewidth\fontdimen8\tenln \@halfwidth .5\@wholewidth}
\def\thicklines{\let\@linefnt\tenlnw \let\@circlefnt\tencircw
  \@wholewidth\fontdimen8\tenlnw \@halfwidth .5\@wholewidth}

\font\tenln    = line10
\font\tenlnw   = linew10
\font\tencirc  = lcircle10    
\font\tencircw = lcirclew10   
\thinlines

\def\@ifstar#1#2{\@ifnextchar *{\def\@tempa*{#1}\@tempa}{#2}}
\def\@ifnextchar#1#2#3{\let\@tempe #1\def\@tempa{#2}\def\@tempb{#3}\futurelet
    \@tempc\@ifnch}
\def\@ifnch{\ifx \@tempc \@sptoken \let\@tempd\@xifnch
      \else \ifx \@tempc \@tempe\let\@tempd\@tempa\else\let\@tempd\@tempb\fi
      \fi \@tempd}

\def\@getcirc#1{\@tempdima #1\relax \advance\@tempdima 2pt\relax
  \@tempcnta\@tempdima
  \@tempdima 4pt\relax \divide\@tempcnta\@tempdima
  \ifnum \@tempcnta > 10\relax \@tempcnta 10\relax\fi
  \ifnum \@tempcnta >\z@ \advance\@tempcnta\m@ne
    \else \@warning{Oval too small}\fi
  \multiply\@tempcnta 4\relax
  \setbox \@tempboxa \hbox{\@circlefnt
  \char \@tempcnta}\@tempdima \wd \@tempboxa}

\def\@put#1#2#3{\raise #2\hbox to \z@{\hskip #1#3\hss}}

\def\circle{\@ifstar{\@dot}{\@circle}}
\def\@circle#1{\begingroup \boxmaxdepth \maxdimen \@tempdimb #1\unitlength
   \ifdim \@tempdimb >15.5pt\relax \@getcirc\@tempdimb
      \@ovro\ht\@tempboxa
     \setbox\@tempboxa\hbox{\@circlefnt
      \advance\@tempcnta\tw@ \char \@tempcnta
      \advance\@tempcnta\m@ne \char \@tempcnta \kern -2\@tempdima
      \advance\@tempcnta\tw@
      \raise \@tempdima \hbox{\char\@tempcnta}\raise \@tempdima
        \box\@tempboxa}\ht\@tempboxa\z@ \dp\@tempboxa\z@
      \@put{-\@ovro}{-\@ovro}{\box\@tempboxa}%
   \else  \@circ\@tempdimb{96}\fi\endgroup}

\def\@dot#1{\@tempdimb #1\unitlength \@circ\@tempdimb{112}}

\def\@circ#1#2{\@tempdima #1\relax \advance\@tempdima .5pt\relax
   \@tempcnta\@tempdima \@tempdima 1pt\relax
   \divide\@tempcnta\@tempdima
   \ifnum\@tempcnta > 15\relax \@tempcnta 15\relax \fi
   \ifnum \@tempcnta >\z@ \advance\@tempcnta\m@ne\fi
   \advance\@tempcnta #2\relax
   \@circlefnt \char\@tempcnta}

\def\@nnil{\@nil}
\def\@empty{}
\def\@fornoop#1\@@#2#3{}

\def\@tfor#1:=#2\do#3{\xdef\@fortmp{#2}\ifx\@fortmp\@empty \else
    \@tforloop#2\@nil\@nil\@@#1{#3}\fi}
\def\@tforloop#1#2\@@#3#4{\def#3{#1}\ifx #3\@nnil
       \let\@nextwhile=\@fornoop \else
      #4\relax\let\@nextwhile=\@tforloop\fi\@nextwhile#2\@@#3{#4}}

\def\@height{height}
\def\@depth{depth}
\def\@width{width}

\def\newbox{\alloc@4\box\chardef\insc@unt}

\def\@makebox[#1]{\leavevmode\@ifnextchar [{\@imakebox[#1]}{\@imakebox[#1][x]}}

\long\def\@imakebox[#1][#2]#3{\hbox to#1{\let\mb@l\hss
\let\mb@r\hss \expandafter\let\csname mb@#2\endcsname\relax
\mb@l #3\mb@r}}

\def\@makepicbox(#1,#2){\leavevmode\@ifnextchar
   [{\@imakepicbox(#1,#2)}{\@imakepicbox(#1,#2)[]}}

\long\def\@imakepicbox(#1,#2)[#3]#4{\vbox to#2\unitlength
   {\let\mb@b\vss \let\mb@l\hss\let\mb@r\hss
    \let\mb@t\vss
    \@tfor\@tempa :=#3\do{\expandafter\let
        \csname mb@\@tempa\endcsname\relax}%
\mb@t\hbox to #1\unitlength{\mb@l #4\mb@r}\mb@b}}

\def\newsavebox#1{\@ifdefinable#1{\newbox#1}}

\def\savebox#1{\@ifnextchar ({\@savepicbox#1}{\@ifnextchar
     [{\@savebox#1}{\sbox#1}}}

\def\sbox#1#2{\setbox#1\hbox{#2}}

\def\@savebox#1[#2]{\@ifnextchar [{\@isavebox#1[#2]}{\@isavebox#1[#2][x]}}

\long\def\@isavebox#1[#2][#3]#4{\setbox#1 \hbox{\@imakebox[#2][#3]{#4}}}

\def\@savepicbox#1(#2,#3){\@ifnextchar
   [{\@isavepicbox#1(#2,#3)}{\@isavepicbox#1(#2,#3)[]}}

\long\def\@isavepicbox#1(#2,#3)[#4]#5{\setbox#1 \hbox{\@imakepicbox
     (#2,#3)[#4]{#5}}}
\long\def\frame#1{\leavevmode
    \hbox{\hskip-\@wholewidth
     \vbox{\vskip-\@wholewidth
            \hrule \@height\@wholewidth
          \hbox{\vrule \@width\@wholewidth #1\vrule \@width\@wholewidth}\hrule
           \@height \@wholewidth\vskip -\@halfwidth}\hskip-\@wholewidth}}

\newdimen\fboxrule
\newdimen\fboxsep

\long\def\fbox#1{\leavevmode\setbox\@tempboxa\hbox{#1}\@tempdima\fboxrule
    \advance\@tempdima \fboxsep \advance\@tempdima \dp\@tempboxa
   \hbox{\lower \@tempdima\hbox
  {\vbox{\hrule \@height \fboxrule
          \hbox{\vrule \@width \fboxrule \hskip\fboxsep
          \vbox{\vskip\fboxsep \box\@tempboxa\vskip\fboxsep}\hskip
                 \fboxsep\vrule \@width \fboxrule}
                 \hrule \@height \fboxrule}}}}

\def\framebox{\@ifnextchar ({\@framepicbox}{\@ifnextchar
     [{\@framebox}{\fbox}}}

\def\@framebox[#1]{\@ifnextchar [{\@iframebox[#1]}{\@iframebox[#1][x]}}

\long\def\@iframebox[#1][#2]#3{\leavevmode
  \savebox\@tempboxa[#1][#2]{\kern\fboxsep #3\kern\fboxsep}\@tempdima\fboxrule
    \advance\@tempdima \fboxsep \advance\@tempdima \dp\@tempboxa
   \hbox{\lower \@tempdima\hbox
  {\vbox{\hrule \@height \fboxrule
          \hbox{\vrule \@width \fboxrule \hskip-\fboxrule
              \vbox{\vskip\fboxsep \box\@tempboxa\vskip\fboxsep}\hskip
                  -\fboxrule\vrule \@width \fboxrule}
                  \hrule \@height \fboxrule}}}}

\def\@framepicbox(#1,#2){\@ifnextchar
   [{\@iframepicbox(#1,#2)}{\@iframepicbox(#1,#2)[]}}

\long\def\@iframepicbox(#1,#2)[#3]#4{\frame{\@imakepicbox(#1,#2)[#3]{#4}}}

\def\parbox{\@ifnextchar [{\@iparbox}{\@iparbox[c]}}

\long\def\@iparbox[#1]#2#3{\leavevmode \@pboxswfalse
   \if #1b\vbox
     \else \if #1t\vtop
              \else \ifmmode \vcenter
                        \else \@pboxswtrue $\vcenter
                     \fi
           \fi
    \fi{\hsize #2\@parboxrestore #3}\if@pboxsw $\fi}

\def\KREIS(#1,#2){\unskip
  \raise#2 \Einheit\hbox to0pt{\hskip#1 \Einheit
          \raise-.5pt\hbox to0pt{\hss\hskip8pt\circle{8}\hss}\hss}}

\def\Line@(#1,#2)#3{\@xarg #1\relax \@yarg #2\relax
\@linelen=#3\Einheit
\ifnum\@xarg =0 \@vline
  \else \ifnum\@yarg =0 \@hline \else \@sline\fi
\fi}

\def\@sline{\ifnum\@xarg< 0 \@negargtrue \@xarg -\@xarg \@yyarg -\@yarg
  \else \@negargfalse \@yyarg \@yarg \fi
\ifnum \@yyarg >0 \@tempcnta\@yyarg \else \@tempcnta -\@yyarg \fi
\ifnum\@tempcnta>6 \@badlinearg\@tempcnta0 \fi
\ifnum\@xarg>6 \@badlinearg\@xarg 1 \fi
\setbox\@linechar\hbox{\@linefnt\@getlinechar(\@xarg,\@yyarg)}%
\ifnum \@yarg >0 \let\@upordown\raise \@clnht\z@
   \else\let\@upordown\lower \@clnht \ht\@linechar\fi
\@clnwd=\wd\@linechar
\if@negarg \hskip -\wd\@linechar \def\@tempa{\hskip -2\wd\@linechar}\else
     \let\@tempa\relax \fi
\@whiledim \@clnwd <\@linelen \do
  {\@upordown\@clnht\copy\@linechar
   \@tempa
   \advance\@clnht \ht\@linechar
   \advance\@clnwd \wd\@linechar}%
\advance\@clnht -\ht\@linechar
\advance\@clnwd -\wd\@linechar
\@tempdima\@linelen\advance\@tempdima -\@clnwd
\@tempdimb\@tempdima\advance\@tempdimb -\wd\@linechar
\if@negarg \hskip -\@tempdimb \else \hskip \@tempdimb \fi
\multiply\@tempdima \@m
\@tempcnta \@tempdima \@tempdima \wd\@linechar \divide\@tempcnta \@tempdima
\@tempdima \ht\@linechar \multiply\@tempdima \@tempcnta
\divide\@tempdima \@m
\advance\@clnht \@tempdima
\ifdim \@linelen <\wd\@linechar
   \hskip \wd\@linechar
  \else\@upordown\@clnht\copy\@linechar\fi}

\def\@hline{\ifnum \@xarg <0 \hskip -\@linelen \fi
\vrule height\Pfadd@cke width \@linelen depth\Pfadd@cke
\ifnum \@xarg <0 \hskip -\@linelen \fi}

\def\@getlinechar(#1,#2){\@tempcnta#1\relax\multiply\@tempcnta 8
\advance\@tempcnta -9 \ifnum #2>0 \advance\@tempcnta #2\relax\else
\advance\@tempcnta -#2\relax\advance\@tempcnta 64 \fi
\char\@tempcnta}

\def\Vektor(#1,#2)#3(#4,#5){\unskip\leavevmode
  \xcoord#4\relax \ycoord#5\relax
      \raise\ycoord \Einheit\hbox to0pt{\hskip\xcoord \Einheit
         \Vector@(#1,#2){#3}\hss}}

\def\Vector@(#1,#2)#3{\@xarg #1\relax \@yarg #2\relax
\@tempcnta \ifnum\@xarg<0 -\@xarg\else\@xarg\fi
\ifnum\@tempcnta<5\relax
\@linelen=#3\Einheit
\ifnum\@xarg =0 \@vvector
  \else \ifnum\@yarg =0 \@hvector \else \@svector\fi
\fi
\else\@badlinearg\fi}

\def\@hvector{\@hline\hbox to 0pt{\@linefnt
\ifnum \@xarg <0 \@getlarrow(1,0)\hss\else
    \hss\@getrarrow(1,0)\fi}}

\def\@vvector{\ifnum \@yarg <0 \@downvector \else \@upvector \fi}

\def\@svector{\@sline
\@tempcnta\@yarg \ifnum\@tempcnta <0 \@tempcnta=-\@tempcnta\fi
\ifnum\@tempcnta <5
  \hskip -\wd\@linechar
  \@upordown\@clnht \hbox{\@linefnt  \if@negarg
  \@getlarrow(\@xarg,\@yyarg) \else \@getrarrow(\@xarg,\@yyarg) \fi}%
\else\@badlinearg\fi}

\def\@upline{\hbox to \z@{\hskip -.5\Pfadd@cke \vrule width \Pfadd@cke
   height \@linelen depth \z@\hss}}

\def\@downline{\hbox to \z@{\hskip -.5\Pfadd@cke \vrule width \Pfadd@cke
   height \z@ depth \@linelen \hss}}

\def\@upvector{\@upline\setbox\@tempboxa\hbox{\@linefnt\char'66}\raise
     \@linelen \hbox to\z@{\lower \ht\@tempboxa\box\@tempboxa\hss}}

\def\@downvector{\@downline\lower \@linelen
      \hbox to \z@{\@linefnt\char'77\hss}}

\def\@getlarrow(#1,#2){\ifnum #2 =\z@ \@tempcnta='33\else
\@tempcnta=#1\relax\multiply\@tempcnta \sixt@@n \advance\@tempcnta
-9 \@tempcntb=#2\relax\multiply\@tempcntb \tw@
\ifnum \@tempcntb >0 \advance\@tempcnta \@tempcntb\relax
\else\advance\@tempcnta -\@tempcntb\advance\@tempcnta 64
\fi\fi\char\@tempcnta}

\def\@getrarrow(#1,#2){\@tempcntb=#2\relax
\ifnum\@tempcntb < 0 \@tempcntb=-\@tempcntb\relax\fi
\ifcase \@tempcntb\relax \@tempcnta='55 \or
\ifnum #1<3 \@tempcnta=#1\relax\multiply\@tempcnta
24 \advance\@tempcnta -6 \else \ifnum #1=3 \@tempcnta=49
\else\@tempcnta=58 \fi\fi\or
\ifnum #1<3 \@tempcnta=#1\relax\multiply\@tempcnta
24 \advance\@tempcnta -3 \else \@tempcnta=51\fi\or
\@tempcnta=#1\relax\multiply\@tempcnta
\sixt@@n \advance\@tempcnta -\tw@ \else
\@tempcnta=#1\relax\multiply\@tempcnta
\sixt@@n \advance\@tempcnta 7 \fi\ifnum #2<0 \advance\@tempcnta 64 \fi
\char\@tempcnta}

\def\Diagonale(#1,#2)#3{\unskip\leavevmode
  \xcoord#1\relax \ycoord#2\relax
      \raise\ycoord \Einheit\hbox to0pt{\hskip\xcoord \Einheit
         \Line@(1,1){#3}\hss}}
\def\AntiDiagonale(#1,#2)#3{\unskip\leavevmode
  \xcoord#1\relax \ycoord#2\relax 
      \raise\ycoord \Einheit\hbox to0pt{\hskip\xcoord \Einheit
         \Line@(1,-1){#3}\hss}}
\def\Pfad(#1,#2),#3\endPfad{\unskip\leavevmode
  \xcoord#1 \ycoord#2 \thicklines\ZeichnePfad#3\endPfad\thinlines}
\def\ZeichnePfad#1{\ifx#1\endPfad\let\next\relax
  \else\let\next\ZeichnePfad
    \ifnum#1=1
      \raise\ycoord \Einheit\hbox to0pt{\hskip\xcoord \Einheit
         \vrule height\Pfadd@cke width1 \Einheit depth\Pfadd@cke\hss}%
      \advance\xcoord by 1
    \else\ifnum#1=2
      \raise\ycoord \Einheit\hbox to0pt{\hskip\xcoord \Einheit
        \hbox{\hskip-\PfadD@cke\vrule height1 \Einheit width\PfadD@cke depth0pt}\hss}%
      \advance\ycoord by 1
    \else\ifnum#1=3
      \raise\ycoord \Einheit\hbox to0pt{\hskip\xcoord \Einheit
         \Line@(1,1){1}\hss}
      \advance\xcoord by 1
      \advance\ycoord by 1
    \else\ifnum#1=4
      \raise\ycoord \Einheit\hbox to0pt{\hskip\xcoord \Einheit
         \Line@(1,-1){1}\hss}
      \advance\xcoord by 1
      \advance\ycoord by -1
    \fi\fi\fi\fi
  \fi\next}
\def\hSSchritt{\leavevmode\raise-.4pt\hbox to0pt{\hss.\hss}\hskip.2\Einheit
  \raise-.4pt\hbox to0pt{\hss.\hss}\hskip.2\Einheit
  \raise-.4pt\hbox to0pt{\hss.\hss}\hskip.2\Einheit
  \raise-.4pt\hbox to0pt{\hss.\hss}\hskip.2\Einheit
  \raise-.4pt\hbox to0pt{\hss.\hss}\hskip.2\Einheit}
\def\vSSchritt{\vbox{\baselineskip.2\Einheit\lineskiplimit0pt
\hbox{.}\hbox{.}\hbox{.}\hbox{.}\hbox{.}}}
\def\DSSchritt{\leavevmode\raise-.4pt\hbox to0pt{%
  \hbox to0pt{\hss.\hss}\hskip.2\Einheit
  \raise.2\Einheit\hbox to0pt{\hss.\hss}\hskip.2\Einheit
  \raise.4\Einheit\hbox to0pt{\hss.\hss}\hskip.2\Einheit
  \raise.6\Einheit\hbox to0pt{\hss.\hss}\hskip.2\Einheit
  \raise.8\Einheit\hbox to0pt{\hss.\hss}\hss}}
\def\dSSchritt{\leavevmode\raise-.4pt\hbox to0pt{%
  \hbox to0pt{\hss.\hss}\hskip.2\Einheit
  \raise-.2\Einheit\hbox to0pt{\hss.\hss}\hskip.2\Einheit
  \raise-.4\Einheit\hbox to0pt{\hss.\hss}\hskip.2\Einheit
  \raise-.6\Einheit\hbox to0pt{\hss.\hss}\hskip.2\Einheit
  \raise-.8\Einheit\hbox to0pt{\hss.\hss}\hss}}
\def\SPfad(#1,#2),#3\endSPfad{\unskip\leavevmode
  \xcoord#1 \ycoord#2 \ZeichneSPfad#3\endSPfad}
\def\ZeichneSPfad#1{\ifx#1\endSPfad\let\next\relax
  \else\let\next\ZeichneSPfad
    \ifnum#1=1
      \raise\ycoord \Einheit\hbox to0pt{\hskip\xcoord \Einheit
         \hSSchritt\hss}%
      \advance\xcoord by 1
    \else\ifnum#1=2
      \raise\ycoord \Einheit\hbox to0pt{\hskip\xcoord \Einheit
        \hbox{\hskip-2pt \vSSchritt}\hss}%
      \advance\ycoord by 1
    \else\ifnum#1=3
      \raise\ycoord \Einheit\hbox to0pt{\hskip\xcoord \Einheit
         \DSSchritt\hss}
      \advance\xcoord by 1
      \advance\ycoord by 1
    \else\ifnum#1=4
      \raise\ycoord \Einheit\hbox to0pt{\hskip\xcoord \Einheit
         \dSSchritt\hss}
      \advance\xcoord by 1
      \advance\ycoord by -1
    \fi\fi\fi\fi
  \fi\next}
\def\Koordinatenachsen(#1,#2){\unskip
 \hbox to0pt{\hskip-.5pt\vrule height#2 \Einheit width.5pt depth1 \Einheit}%
 \hbox to0pt{\hskip-1 \Einheit \xcoord#1 \advance\xcoord by1
    \vrule height0.25pt width\xcoord \Einheit depth0.25pt\hss}}
\def\Koordinatenachsen(#1,#2)(#3,#4){\unskip
 \hbox to0pt{\hskip-.5pt \ycoord-#4 \advance\ycoord by1
    \vrule height#2 \Einheit width.5pt depth\ycoord \Einheit}%
 \hbox to0pt{\hskip-1 \Einheit \hskip#3\Einheit
    \xcoord#1 \advance\xcoord by1 \advance\xcoord by-#3
    \vrule height0.25pt width\xcoord \Einheit depth0.25pt\hss}}
\def\Gitter(#1,#2){\unskip \xcoord0 \ycoord0 \leavevmode
  \LOOP\ifnum\ycoord<#2
    \loop\ifnum\xcoord<#1
      \raise\ycoord \Einheit\hbox to0pt{\hskip\xcoord \Einheit\Punkt\hss}%
      \advance\xcoord by1
    \repeat
    \xcoord0
    \advance\ycoord by1
  \REPEAT}
\def\Gitter(#1,#2)(#3,#4){\unskip \xcoord#3 \ycoord#4 \leavevmode
  \LOOP\ifnum\ycoord<#2
    \loop\ifnum\xcoord<#1
      \raise\ycoord \Einheit\hbox to0pt{\hskip\xcoord \Einheit\Punkt\hss}%
      \advance\xcoord by1
    \repeat
    \xcoord#3
    \advance\ycoord by1
  \REPEAT}
\def\Label#1#2(#3,#4){\unskip \xdim#3 \Einheit \ydim#4 \Einheit
  \def\lo{\advance\xdim by-.5 \Einheit \advance\ydim by.5 \Einheit}%
  \def\llo{\advance\xdim by-.25cm \advance\ydim by.5 \Einheit}%
  \def\loo{\advance\xdim by-.5 \Einheit \advance\ydim by.25cm}%
  \def\o{\advance\ydim by.25cm}%
  \def\ro{\advance\xdim by.5 \Einheit \advance\ydim by.5 \Einheit}%
  \def\rro{\advance\xdim by.25cm \advance\ydim by.5 \Einheit}%
  \def\roo{\advance\xdim by.5 \Einheit \advance\ydim by.25cm}%
  \def\l{\advance\xdim by-.30cm}%
  \def\r{\advance\xdim by.30cm}%
  \def\lu{\advance\xdim by-.5 \Einheit \advance\ydim by-.6 \Einheit}%
  \def\llu{\advance\xdim by-.25cm \advance\ydim by-.6 \Einheit}%
  \def\luu{\advance\xdim by-.5 \Einheit \advance\ydim by-.30cm}%
  \def\u{\advance\ydim by-.30cm}%
  \def\ru{\advance\xdim by.5 \Einheit \advance\ydim by-.6 \Einheit}%
  \def\rru{\advance\xdim by.25cm \advance\ydim by-.6 \Einheit}%
  \def\ruu{\advance\xdim by.5 \Einheit \advance\ydim by-.30cm}%
  #1\raise\ydim\hbox to0pt{\hskip\xdim
     \vbox to0pt{\vss\hbox to0pt{\hss$#2$\hss}\vss}\hss}%
}

\font@\twelverm=cmr10 scaled\magstep1
\font@\twelveit=cmti10 scaled\magstep1
\font@\twelvebf=cmbx10 scaled\magstep1
\font@\twelvei=cmmi10 scaled\magstep1
\font@\twelvesy=cmsy10 scaled\magstep1
\font@\twelveex=cmex10 scaled\magstep1

\newtoks\twelvepoint@
\def\twelvepoint{\normalbaselineskip15\p@
 \abovedisplayskip15\p@ plus3.6\p@ minus10.8\p@
 \belowdisplayskip\abovedisplayskip
 \abovedisplayshortskip\z@ plus3.6\p@
 \belowdisplayshortskip8.4\p@ plus3.6\p@ minus4.8\p@
 \textonlyfont@\rm\twelverm \textonlyfont@\it\twelveit
 \textonlyfont@\sl\twelvesl \textonlyfont@\bf\twelvebf
 \textonlyfont@\smc\twelvesmc \textonlyfont@\tt\twelvett
%
 \ifsyntax@ \def\big##1{{\hbox{$\left##1\right.$}}}%
  \let\Big\big \let\bigg\big \let\Bigg\big
 \else
  \textfont\z@=\twelverm  \scriptfont\z@=\tenrm  \scriptscriptfont\z@=\sevenrm
  \textfont\@ne=\twelvei  \scriptfont\@ne=\teni  \scriptscriptfont\@ne=\seveni
  \textfont\tw@=\twelvesy \scriptfont\tw@=\tensy \scriptscriptfont\tw@=\sevensy
  \textfont\thr@@=\twelveex \scriptfont\thr@@=\tenex
        \scriptscriptfont\thr@@=\tenex
  \textfont\itfam=\twelveit \scriptfont\itfam=\tenit
        \scriptscriptfont\itfam=\tenit
  \textfont\bffam=\twelvebf \scriptfont\bffam=\tenbf
        \scriptscriptfont\bffam=\sevenbf
  \setbox\strutbox\hbox{\vrule height10.2\p@ depth4.2\p@ width\z@}%
  \setbox\strutbox@\hbox{\lower.6\normallineskiplimit\vbox{%
        \kern-\normallineskiplimit\copy\strutbox}}%
 \setbox\z@\vbox{\hbox{$($}\kern\z@}\bigsize@=1.4\ht\z@
 \fi
 \normalbaselines\rm\ex@.2326ex\jot3.6\ex@\the\twelvepoint@}

\font@\fourteenrm=cmr10 scaled\magstep2
\font@\fourteenit=cmti10 scaled\magstep2
\font@\fourteensl=cmsl10 scaled\magstep2
\font@\fourteensmc=cmcsc10 scaled\magstep2
\font@\fourteentt=cmtt10 scaled\magstep2
\font@\fourteenbf=cmbx10 scaled\magstep2
\font@\fourteeni=cmmi10 scaled\magstep2
\font@\fourteensy=cmsy10 scaled\magstep2
\font@\fourteenex=cmex10 scaled\magstep2
\font@\fourteenmsa=msam10 scaled\magstep2
\font@\fourteeneufm=eufm10 scaled\magstep2
\font@\fourteenmsb=msbm10 scaled\magstep2
\newtoks\fourteenpoint@
\def\fourteenpoint{\normalbaselineskip15\p@
 \abovedisplayskip18\p@ plus4.3\p@ minus12.9\p@
 \belowdisplayskip\abovedisplayskip
 \abovedisplayshortskip\z@ plus4.3\p@
 \belowdisplayshortskip10.1\p@ plus4.3\p@ minus5.8\p@
 \textonlyfont@\rm\fourteenrm \textonlyfont@\it\fourteenit
 \textonlyfont@\sl\fourteensl \textonlyfont@\bf\fourteenbf
 \textonlyfont@\smc\fourteensmc \textonlyfont@\tt\fourteentt
%
 \ifsyntax@ \def\big##1{{\hbox{$\left##1\right.$}}}%
  \let\Big\big \let\bigg\big \let\Bigg\big
 \else
  \textfont\z@=\fourteenrm  \scriptfont\z@=\twelverm  \scriptscriptfont\z@=\tenrm
  \textfont\@ne=\fourteeni  \scriptfont\@ne=\twelvei  \scriptscriptfont\@ne=\teni
  \textfont\tw@=\fourteensy \scriptfont\tw@=\twelvesy \scriptscriptfont\tw@=\tensy
  \textfont\thr@@=\fourteenex \scriptfont\thr@@=\twelveex
        \scriptscriptfont\thr@@=\twelveex
  \textfont\itfam=\fourteenit \scriptfont\itfam=\twelveit
        \scriptscriptfont\itfam=\twelveit
  \textfont\bffam=\fourteenbf \scriptfont\bffam=\twelvebf
        \scriptscriptfont\bffam=\tenbf
  \setbox\strutbox\hbox{\vrule height12.2\p@ depth5\p@ width\z@}%
  \setbox\strutbox@\hbox{\lower.72\normallineskiplimit\vbox{%
        \kern-\normallineskiplimit\copy\strutbox}}%
 \setbox\z@\vbox{\hbox{$($}\kern\z@}\bigsize@=1.7\ht\z@
 \fi
 \normalbaselines\rm\ex@.2326ex\jot4.3\ex@\the\fourteenpoint@}

\catcode`\@=13
\def\[{\left[}
\def\]{\right]}
%

\topmatter
\title 
A factorization theorem for lozenge tilings of a hexagon with triangular holes
\endtitle
\author M.~Ciucu$^\dagger$ and C.~Krattenthaler$^\ddagger$
\endauthor
\affil Department of Mathematics,
Indiana University,\\
Bloomington, IN 47405-5701, USA\\\vskip6pt
Fakult\"at f\"ur Mathematik der Universit\"at Wien,\\
Oskar-Morgenstern-Platz 1, A-1090 Vienna, Austria.\\
WWW: {\tt http://www.mat.univie.ac.at/\lower0.5ex\hbox{\~{}}kratt}
\endaffil
\address Department of Mathematics, Indiana University, Bloomington, IN 47405-5701,
USA
\endaddress
\address Fakult\"at f\"ur Mathematik der Universit\"at Wien,
Oskar-Morgenstern-Platz 1, A-1090 Vienna, Austria.
\endaddress
\thanks{$^\dagger$Research partially supported by NSF grants DMS-1101670
and DMS-1501052.}\endthanks
\thanks{$^\ddagger$Research partially supported by the Austrian
Science Foundation FWF, grants Z130-N13 and F50-N15,
the latter in the framework of the Special Research Program
``Algorithmic and Enumerative Combinatorics."}\endthanks
\abstract
In this paper we present a combinatorial generalization of the fact
that the number of plane partitions that fit in a $2a\times b\times b$
box is equal to the number of such plane partitions that are
symmetric, times the number of such plane partitions for which the
transpose is the same as the complement. We use the equivalent
phrasing of this identity in terms of symmetry classes of lozenge
tilings of a hexagon on the triangular lattice. Our generalization
consists of allowing the hexagon have certain symmetrically placed
holes along its horizontal symmetry axis. The special case when there
are no holes can be viewed as a new, simpler proof of the enumeration
of symmetric plane partitions.  
\endabstract
\endtopmatter
\document

\leftheadtext{M. Ciucu and C. Krattenthaler}

\rightheadtext{A factorization theorem for lozenge tilings}

\head 1. Introduction\endhead

The enumeration of the ten symmetry classes of plane partitions that
fit in a given box --- equivalently, symmetry classes of lozenge
tilings of hexagons on the triangular lattice --- forms a classical
chapter of enumerative combinatorics
(see \cite{\StanSymmClass},\cite{\BresAO},\cite{\AndrAW},%
\cite{\KupNnine},\cite{\StemAG},\cite{\KaKZAC}). Explicit
product formulas exist for all symmetry classes, which makes it
possible to find relations between them. One such striking relation is
that   
$$
\M(H_{2a,b,b})=\M_{-}(H_{2a,b,b})\M_{|}(H_{2a,b,b}),\tag\eaa
$$
where $H_{2a,b,b}$ is the hexagon of side-lengths $2a$, $b$, $b$,
$2a$, $b$, $b$ (clockwise from the western\footnote{ Throughout this
paper, we consider the triangular lattice drawn in the plane so that
one of the families of lattice lines is vertical.} side), $\M(R)$
denotes the number of lozenge tilings\footnote{ A {\it lozenge} is the
union of two unit triangles sharing an edge; a {\it lozenge tiling} of
a lattice region $R$ is a covering of $R$ by lozenges, with no holes or
overlaps.} of the lattice region $R$, and $\M_{-}(R)$ (resp.,
$\M_{|}(R)$) is the number of lozenge tilings of $R$ that are
invariant under reflection across the horizontal (resp., vertical)
symmetry axis of $R$ (provided $R$ possesses such symmetries). 

Equation (\eaa) is an immediate consequence of the explicit formulas
enumerating the symmetry classes of plane partitions (in the notation
of \cite{\StanSymmClass}, $\M(H_{2a,b,b})=N_1(2a,b,b)$ is the number
of plane partitions fitting in a $2a\times b\times b$ 
box, $\M_{-}(H_{2a,b,b})=N_6(2a,b,b)$ is the number of {\it
transpose-complementary} plane partitions, and\linebreak
$\M_{|}(H_{2a,b,b})=N_2(2a,b,b)$ is the number of {\it symmetric}
plane partitions, fitting in the same box). However, the
simplicity of (\eaa) raises two natural questions: How can one see
directly (without explicitly evaluating the terms) that the equation
holds? And how can one generalize it? 

We presented a generalization in terms of Schur functions
in \cite{\CiKrAF}, which gave an algebraic reason for why equation
(\eaa) holds. In this paper we present a combinatorial generalization,
in terms of hexagons that are allowed to have certain symmetrically
placed holes along their horizontal symmetry axis. The special case
when there are no holes can be regarded as a new proof of the
enumeration of symmetric plane partitions (first proved by
Andrews \cite{\AndrAK}), as it follows, via our factorization result,
from the base case (due to MacMahon \cite{\MacMAA}) and the 
transpose-complementary case (due to Proctor \cite{\Proc}),
as we explain in Section~6. Our results are
described in the next section. 

There are several other simple equations relating the symmetry classes
of plane partitions which can be proved directly (see
e.g\. \cite{\CiucNnineNten}\cite{\CiKrfiveandahalf}\cite{\KupExploration}). 
There
is still no unified proof available for all ten symmetry classes, but
new direct ways of relating them to one another may help achieving
this goal.

\head 2. The factorization theorem\endhead

\topinsert
\centerline{\mypic{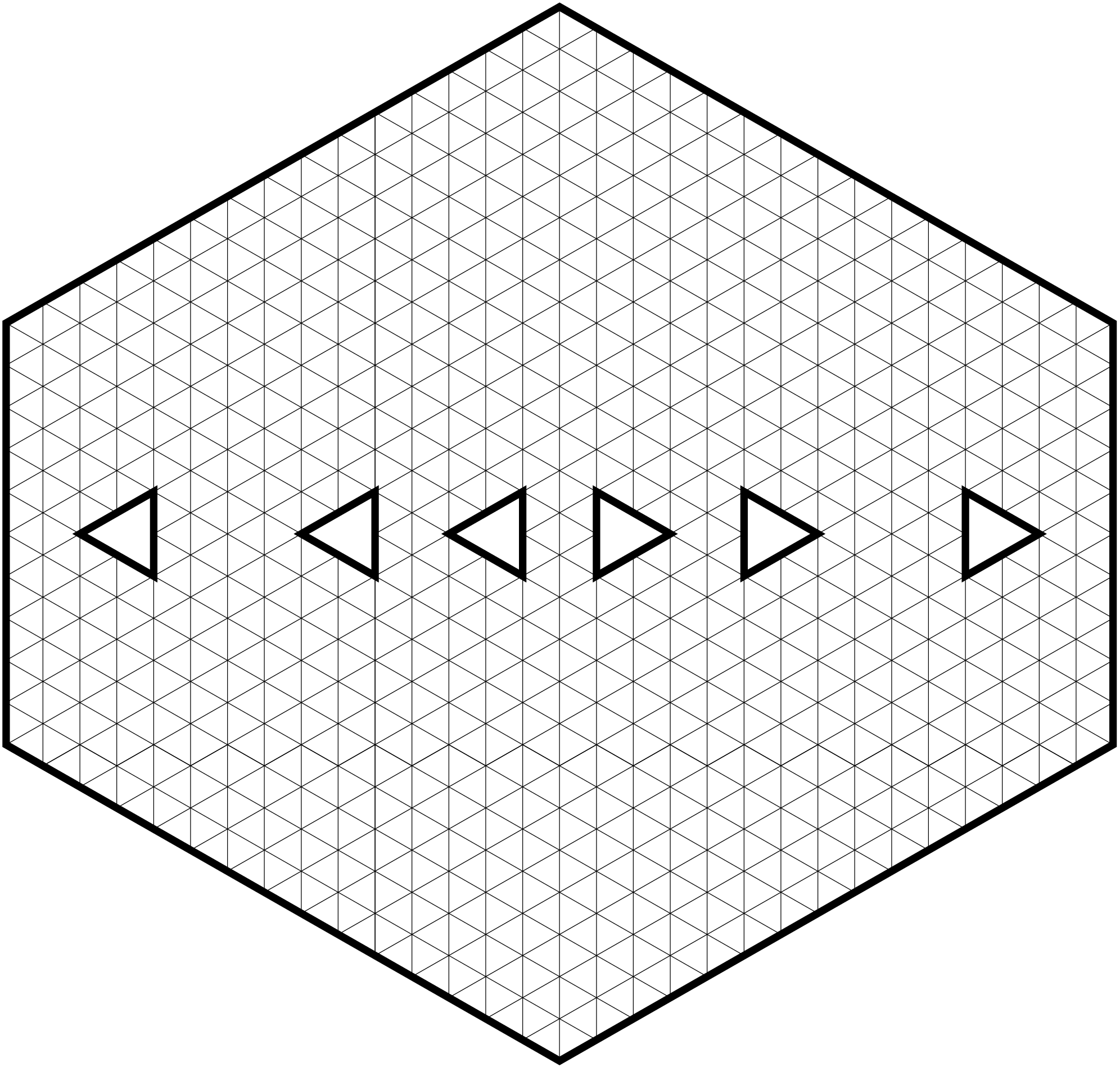}}
\centerline{{\smc Figure~\fba. {\rm The hexagon with holes $H_{15,10}(2,5,7)$.}}}
\endinsert

Let $n,m,l$ be positive integers, and let $k_1,k_2,\dotsc,k_l$ be
positive integers with $k_1<k_2<\dots<k_l\le n/2$.
Denote by $H_{n,2m}(k_1,k_2,\dots,k_l)$ the region obtained from the
hexagon with side lengths $n,2m,n,n,2m,n$ (where the sides of length
$2m$ are vertical) by removing the following $2l$ triangles of side
length two from along its horizontal symmetry axis: $l$ left-pointing
such triangles, with vertical sides at distances
$2k_1,2k_2,\dots,2k_l$ from the left side of the hexagon (in units
equal to $\sqrt{3}$ times the lattice spacing), and their reflections
across the vertical symmetry axis of the hexagon. 
Figure~\fba\ shows the region $H_{15,10}(2,5,7)$.

The factorization theorem which we prove in this paper is the
following generalization of (\eaa).

\pagebreak

\proclaim{Theorem \tba} For all positive integers $n,m,l$ and
positive integers $k_1,k_2,\dots k_l$ with 
$k_1<k_2<\dots<k_l\le n/2$, we have 
$$
\M\left(H_{n,2m}(k_1,k_2,\dots,k_l)\right)
=
\M_{-}\left(H_{n,2m}(k_1,k_2,\dots,k_l)\right)
\M_{|}\left(H_{n,2m}(k_1,k_2,\dots,k_l)\right).
\tag\eba
$$

\endproclaim

There are $n-2l$ lozenge positions (i.e., positions that might be
occupied by a lozenge in some tiling) along the horizontal symmetry
axis of $H_{n,2m}(k_1,k_2,\dots,k_l)$ which can accommodate horizontal
lozenges (these are shaded in the picture on the left in Figure~{\fbb}). 
Clearly, a necessary condition for a tiling to be symmetric
about the horizontal symmetry axis is that all these $n-2l$ positions
are occupied by horizontal lozenges. 

It follows that, if we denote by  $H^+_{n,2m}(k_1,k_2,\dots,k_l)$ the
portion of the region $H_{n,2m}(k_1,k_2,\dots,k_l)$ that is above the
zig-zag line which starts at the center of the left side of the
hexagon and proceeds just above its horizontal symmetry axis until it
reaches the center of the right 
side (for $H_{15,10}(2,5,7)$ this construction is pictured on the left
in Figure~{\fbb}), then we have 
$$
\M_{-}\left(H_{n,2m}(k_1,k_2,\dots,k_l)\right)
=
\M\left(H^+_{n,2m}(k_1,k_2,\dots,k_l)\right).\tag\ebb
$$

Note also that if $F_{n,2m}(k_1,k_2,\dots,k_l)$ is the region
consisting of the left half of $H_{n,2m}(k_1,k_2,\dots,k_l)$, with the
portion of its boundary that is along the vertical symmetry axis of
$H_{n,2m}(k_1,k_2,\dots,k_l)$ taken to be free (i.e., when considering
lozenge tilings of this region, lozenges are allowed to protrude out
halfway across this part of the boundary; $F_{15,10}(2,5,7)$ is
illustrated on the right in Figure~{\fbb}), then we have 
$$
\M_{|}\left(H_{n,2m}(k_1,k_2,\dots,k_l)\right)
=
\M_f\left(F_{n,2m}(k_1,k_2,\dots,k_l)\right),\tag\ebc
$$
where, for a lattice region $R$ with free boundary conditions along
some portion of its boundary, $\M_f(R)$ denotes the number of lozenge
tilings of $R$ in which lozenges are allowed to protrude out halfway
across the free part of the boundary. 

In view of (\ebb) and (\ebc), the statement of Theorem~{\tba} is
equivalent to the equality 
$$
\multline
\M\big(H_{n,2m}(k_1,k_2,\dots,k_l)\big)\\
=\M\big(H^+_{n,2m}(k_1,k_2,\dots,k_l)\big)
\,\M_f\big(F_{n,2m}(k_1,k_2,\dots,k_l)\big).
\endmultline
\tag\AA
$$


On the other hand, if $H^-_{n,2m}(k_1,k_2,\dotsc,k_l)$ is the region
below the zig-zag line that defined $H^+_{n,2m}(k_1,k_2,\dotsc,k_l)$,
with the extra specification that the $n-2l$ lozenge positions fitting
in the ``folds'' of the zig-zag are weighted by $1/2$ (the shaded
positions on the left in Figure~{\fbb}), then the factorization
theorem of \cite{\CiucAB, Theorem~1.2} implies that 
$$
\multline
\M\big(H_{n,2m}(k_1,k_2,\dots,k_l)\big)\\
=2^{n-2l}\,
\M\big(H^+_{n,2m}(k_1,k_2,\dots,k_l)\big)
\,\M^*\big(H^-_{n,2m}(k_1,k_2,\dots,k_l)\big).
\endmultline
\tag\AB
$$
Here, $\M^*\big(H^-_{n,2m}(k_1,k_2,\dots,k_l)\big)$ denotes the
{\it weighted} count of the lozenge tilings of
$H^-_{n,2m}(k_1,k_2,\dots,k_l)$, in which each lozenge occupying one
of the special $n-2l$ tile positions has weight $1/2$, all other
lozenges have weight $1$, and the weight of a tiling is the product of
the weights of its tiles. 

\topinsert
\twolinetwo{\mypic{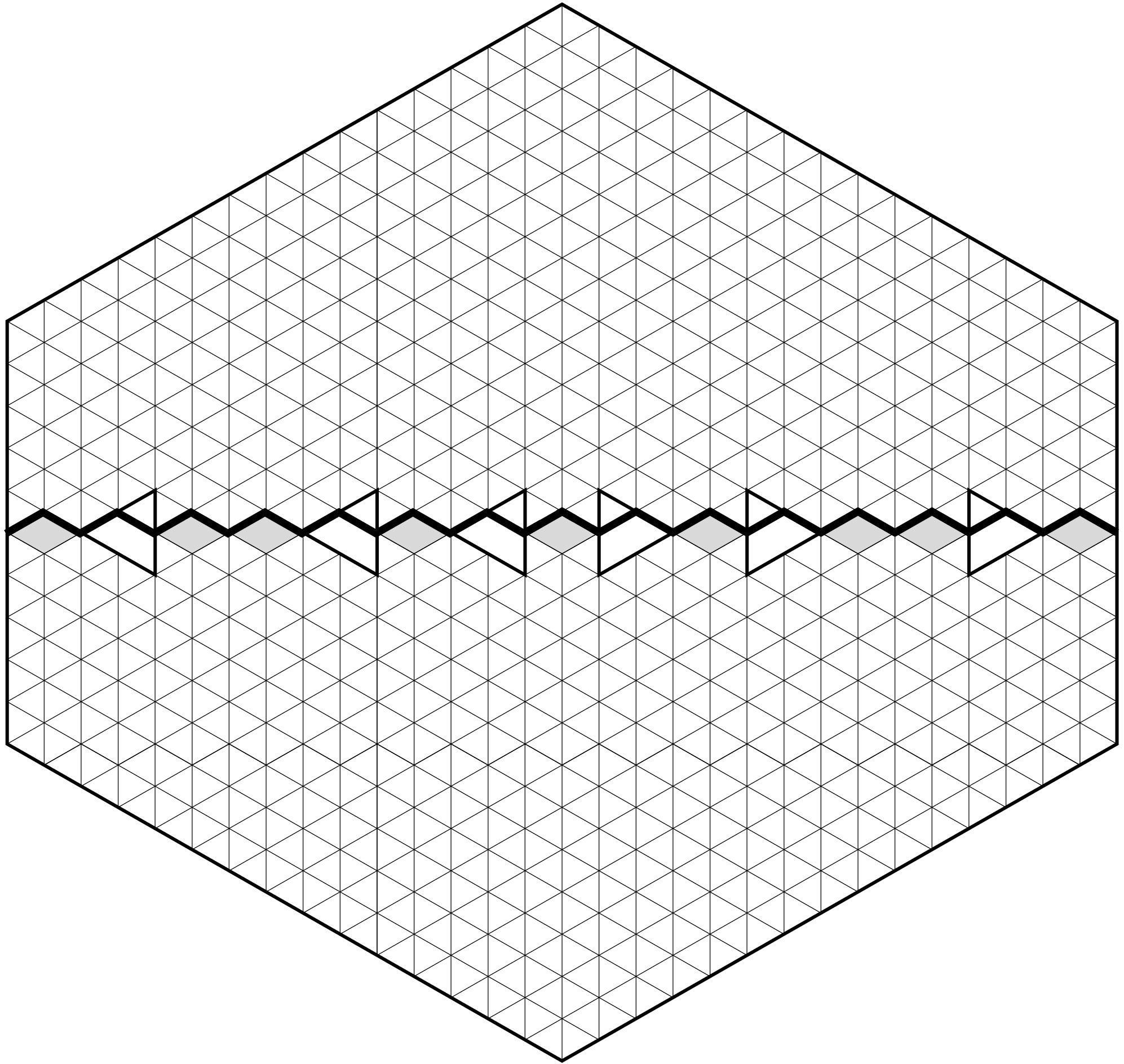}}{\mypic{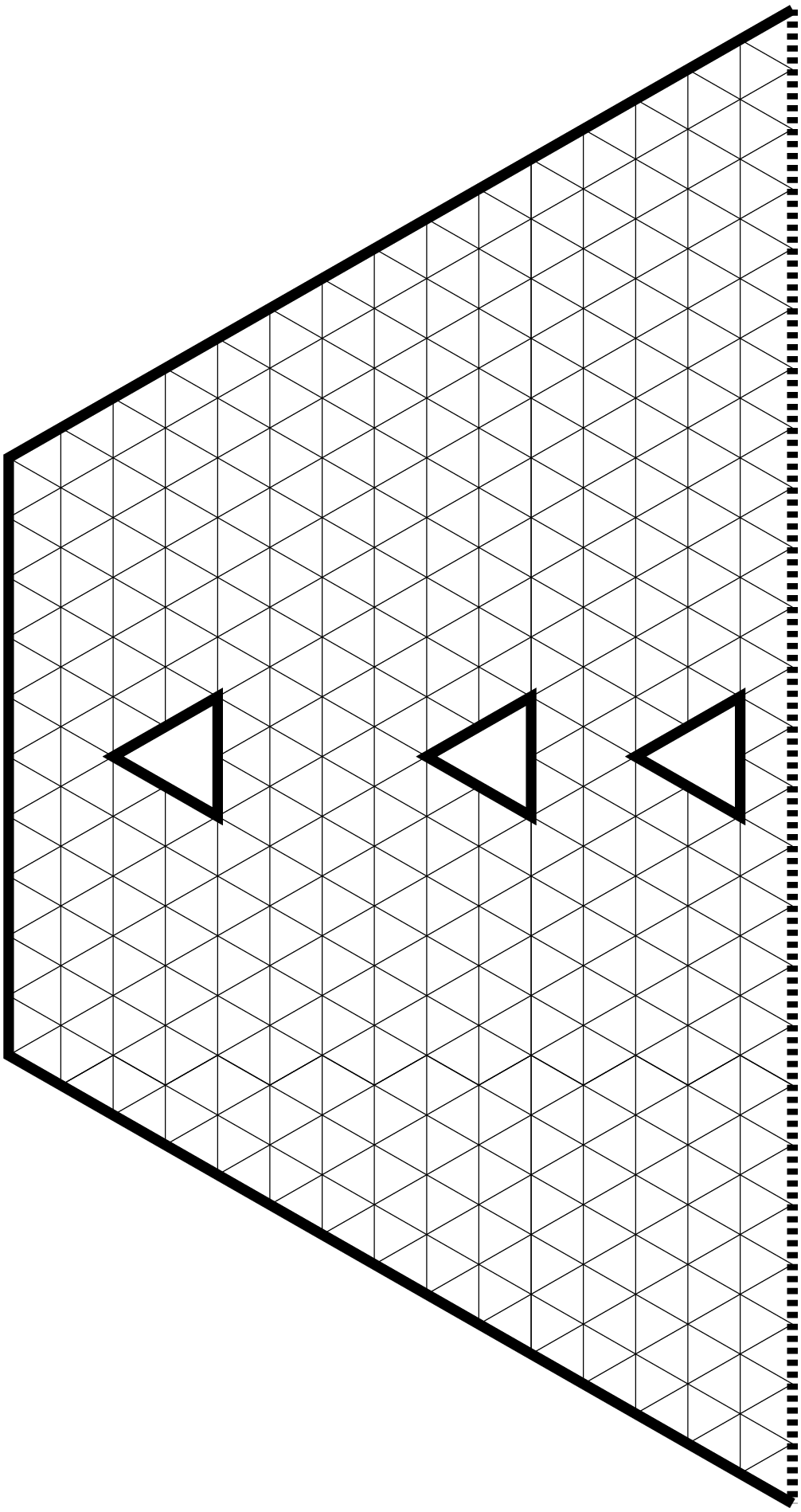}}
\medskip
\centerline{\smc Figure~\fbb. \rm The regions {\rm $H^+_{15,10}(2,5,7)$ and
$H^-_{15,10}(2,5,7)$ (left); $F_{15,10}(2,5,7)$ (right).}} 
\endinsert

Comparison of (\AB) with (\AA) shows that, in order to prove
Theorem~{\tba}, it suffices to establish the relation
$$
\M_f\big(F_{n,2m}(k_1,k_2,\dots,k_l)\big)
=2^{n-2l}\,
\M^*\big(H^-_{n,2m}(k_1,k_2,\dots,k_l)\big).
\tag\AC
$$
This is what we do in the next two sections.

\head 3. 
Non-intersecting lattice paths, Pfaffians, and determinants\endhead

Recall that a family of lattice paths is called {\it non-intersecting}
if no two paths in the family share a vertex. 
In this section, we interpret the perfect matching counts on each side
of (\AC) as the (weighted) 
number of a certain family of non-intersecting lattice paths.  
Using well-known 
determinantal and Pfaffian formulas for the number of non-intersecting
lattice paths, this enables us to express the left-hand side in terms
of a Pfaffian and the right-hand side in terms of a determinant.

We start with the left-hand side. If we apply the standard translation
of lozenge tilings to families of non-intersecting lattice paths
(see e.g\. \cite{\DT} and \cite{\CiucNnineNten}), then we obtain that
$\M_f\left(F_{n,2m}(k_1,k_2,\dots,k_l)\right)$ is equal to the number
of all families 

$$(P_{-m+1},P_{-m+2},\dots,P_{m},P_{1^-},P_{2^-},\dots,P_{l^-},
P_{1^+},P_{2^+},\dots,P_{l^+})$$
of non-intersecting lattice paths,
where for $s\in\{-m+1,-m+2,\dots,m\}$, $P_s$ starts at 
$A_s:=(s,-s+1)$,
while, for $t\in\{1,2,\dots,l\}$, $P_{t^-}$ starts at $B_{t^-}:=(k_t,k_t+1)$
and $P_{t^+}$ starts at $B_{t^+}:=(k_t+1,k_t)$, and all paths end somewhere
on the line $x+y=n+1$. See Figure~\FC\ for the family of
non-intersecting lattice paths which corresponds to the tiling
in Figure~{\FB} (which also indicates, by dotted lines, the
non-intersecting paths of lozenges that encode that tiling). 

\topinsert
\vskip10pt
\vbox{\noindent
\centertexdraw{
\drawdim truecm \setunitscale.7
\linewd.15
\move(0 4)
\RhombusB \RhombusB \RhombusB 
\RhombusA \RhombusB \RhombusA \RhombusB \RhombusA \RhombusA \RhombusB 
\RhombusA \RhombusB \RhombusB \RhombusA \RhombusA \RhombusA
\move(0.866025 4.5)
\RhombusB \RhombusB \RhombusB \RhombusA \RhombusB
\move(1.73205 5)
\RhombusB \RhombusB \RhombusB \RhombusA \RhombusB
\move(2.59808 5.5)
\RhombusB \RhombusB \RhombusA \RhombusB \RhombusB \RhombusA \RhombusB \RhombusB
\RhombusA \RhombusB \RhombusA \RhombusB \RhombusA 
\move(3.4641 6)
\RhombusB \RhombusB \RhombusA \RhombusB \RhombusB \RhombusA \RhombusA \RhombusB
\RhombusB \RhombusB \RhombusA 
\move(4.33013 6.5)
\RhombusB \RhombusB \RhombusA \RhombusB \RhombusA \RhombusB
\RhombusA
\move(5.19615 7)
\RhombusB \RhombusA \RhombusB \RhombusB \RhombusA
\move(6.06218 7.5)
\RhombusB \RhombusA
\move(0 1)
\RhombusC 
\move(0 0)
\RhombusC \RhombusC 
\move(0 -1)
\RhombusC \RhombusC \RhombusC \RhombusC 
\move(0 -2)
\RhombusC \RhombusC \RhombusC \RhombusC \RhombusC 
\move(0 -3)
\RhombusC \RhombusC \RhombusC \RhombusC \RhombusC 
\move(3.4641 1)
\RhombusC 
\move(3.4641 0)
\RhombusC 
\move(2.59808 -.5)
\RhombusC \RhombusC \RhombusC 
\move(4.33013 -2.5)
\RhombusC \RhombusC 
\move(5.19615 -4)
\RhombusC \RhombusC 
\move(5.19615 -5)
\RhombusC \RhombusC 
\move(5.19615 1)
\RhombusC 
\linewd.08
\move(0 -3.5)
\hdSchritt \hdSchritt \hdSchritt \hdSchritt \hdSchritt \odaSchritt
\hdSchritt \hdSchritt 
\move(0 -2.5)
\hdSchritt \hdSchritt \hdSchritt \hdSchritt \hdSchritt \odaSchritt
\hdSchritt \hdSchritt 
\move(0 -1.5)
\hdSchritt \hdSchritt \hdSchritt \hdSchritt \odaSchritt \hdSchritt 
\hdSchritt \odaSchritt 
\move(0 -.5)
\hdSchritt \hdSchritt \odaSchritt \hdSchritt \hdSchritt \hdSchritt 
\odaSchritt \odaSchritt 
\move(0 .5)
\hdSchritt \odaSchritt \odaSchritt \odaSchritt \odaSchritt \odaSchritt 
\hdSchritt \odaSchritt 
\move(0 1.5)
\odaSchritt \odaSchritt \odaSchritt \hdSchritt \odaSchritt \odaSchritt 
\odaSchritt \odaSchritt 
\move(0 2.5)
\odaSchritt \odaSchritt \odaSchritt \odaSchritt \odaSchritt \odaSchritt 
\hdSchritt \odaSchritt 
\move(0 3.5)
\odaSchritt \odaSchritt \odaSchritt \odaSchritt \odaSchritt \odaSchritt 
\odaSchritt \odaSchritt 
\move(3.4641 -.5)
\hdSchritt \odaSchritt \hdSchritt \odaSchritt 
\move(3.4641 .5)
\hdSchritt \odaSchritt \hdSchritt \odaSchritt 
\ringerl(0 -3.5)
\ringerl(0 -2.5)
\ringerl(0 -1.5)
\ringerl(0 -.5)
\ringerl(0 .5)
\ringerl(0 1.5)
\ringerl(0 2.5)
\ringerl(0 3.5)
\ringerl(3.4641 .5)
\ringerl(3.4641 -.5)
\ringerl(6.9282 -6.5)
\ringerl(6.9282 -5.5)
\ringerl(6.9282 -3.5)
\ringerl(6.9282 -1.5)
\ringerl(6.9282 -.5)
\ringerl(6.9282 .5)
\ringerl(6.9282 2.5)
\ringerl(6.9282 4.5)
\ringerl(6.9282 5.5)
\ringerl(6.9282 7.5)
\linewd.05
}
\centerline{\eightpoint A lozenge tiling of the region
$F_{8,8}(2)$; the right boundary 
is free. The}
\centerline{\eightpoint dotted lines mark paths of lozenges. 
They determine the tiling uniquely.}
\vskip8pt
\centerline{\smc Figure \FB}
}
\vskip10pt
\endinsert

\topinsert
\vbox{\noindent
$$
\Gitter(12,13)(-4,-4)
\Koordinatenachsen(12,13)(-4,-4)
\Pfad(-3,4),22222222\endPfad
\Pfad(-2,3),22222212\endPfad
\Pfad(-1,2),22212222\endPfad
\Pfad(0,1),12222212\endPfad
\Pfad(1,0),11211122\endPfad
\Pfad(2,-1),11112112\endPfad
\Pfad(3,-2),11111211\endPfad
\Pfad(4,-3),11111211\endPfad
\Pfad(2,3),1212\endPfad
\Pfad(3,2),1212\endPfad
\DickPunkt(-3,4)
\DickPunkt(-2,3)
\DickPunkt(-1,2)
\DickPunkt(0,1)
\DickPunkt(1,0)
\DickPunkt(2,-1)
\DickPunkt(3,-2)
\DickPunkt(4,-3)
\DickPunkt(2,3)
\DickPunkt(3,2)
\DickPunkt(11,-2)
\DickPunkt(10,-1)
\DickPunkt(8,1)
\DickPunkt(6,3)
\DickPunkt(5,4)
\DickPunkt(4,5)
\DickPunkt(2,7)
\DickPunkt(0,9)
\DickPunkt(-1,10)
\DickPunkt(-3,12)
\Kreis(0,0)
\Label\u{A_{-3}}(-3,4)
\Label\u{A_{-2}}(-2,3)
\Label\u{A_{-1}}(-1,2)
\Label\u{A_{0}}(0,1)
\Label\u{A_1}(1,0)
\Label\u{A_2}(2,-1)
\Label\u{A_3}(3,-2)
\Label\u{A_4}(4,-3)
\Label\o{B_{1^-}}(2,3)
\Label\o{B_{1^+}}(3,2)
\Label\ro{P_{1^-}}(3,4)
\Label\ro{P_{1^+}}(4,3)
\Label\ro{P_4}(9,-4)
\Label\ro{P_3}(7,-2)
\Label\ro{P_2}(5,-1)
\Label\ro{P_1}(2,0)
\Label\o{P_{0}}(2,5)
\Label\u{P_{-1}}(-1,6)
\Label\u{P_{-2}}(-2,10)
\Label\o{P_{-3}}(-4,7)
\Label\r{x}(12,0)
\Label\ru{y}(0,13)
\hskip3cm
$$
\centerline{\eightpoint The paths of lozenges of Figure {\FB} drawn as
non-intersecting lattice paths on $\Z^2$}
\vskip8pt
\centerline{\smc Figure \FC}
}
\vskip10pt
\endinsert

By a slight extension of a theorem due to
Okada \cite{\OkadAA, Theorem~3} and Stembridge \cite{\StemAE, Theorem~3.1},
the number of the above families of non-intersecting lattice paths
can be expressed in terms of a Pfaffian.
The reader should recall that the {\it Pfaffian} of a skew-symmetric
$2n\times 2n$ matrix $A$ can be defined by (see e.g.\ 
\cite{\StemAE, p.~102})
$$\Pf A:=
\sum _{\pi\in\Cal M[1,\dots, 2n]} ^{}\sgn\pi 
\underset i,j\text{ matched in }\pi\to{\prod _{i<j} ^{}}
A_{i,j}\ ,
\tag\AD$$
where $\Cal M[1,2,\dots, 2n]$ denotes the set of all perfect matchings
($1$-factors) of (the complete graph on) $\{1,2,\dots,2n\}$,
$\text{cr}(\pi)$ is the number of crossings\footnote{ A crossing is a
quadruple $i<j<k<l$ such that, under $\pi$, $i$ is paired with $k$,
and $j$ is paired with $l$.} of the perfect matching $\pi$, and
$\sgn\pi=(-1)^{\text{cr}(\pi)}$. 
It is a well-known fact (see e.g\. \cite{\StemAE, Prop.~2.2}) that
$$
(\Pf A)^2=\det A.
\tag\AE$$

\proclaim{Theorem \TC\ \smc (Okada, Stembridge)}
Let $\{u_1,u_2,\dots,u_{p}\}$ and
$I=\{I_1,I_2,\dots\}$ be finite sets of lattice points in the integer lattice
$\Z^2$, with $p$ even. Let $\frak S_p$ be the symmetric group on
$\{1,2,\dots,p\}$, set 
$\bold u_\pi=(u_{\pi(1)},u_{\pi(2)},\dots,u_{\pi(p)})$,
and denote by $\Cal P^{\text{nonint}}(\bold u_\pi\to I)$ 
the number of families $(P_1,P_2,\dots,P_{p})$ of
non-intersecting lattice paths consisting of unit horizontal and
vertical steps in the positive direction, with $P_k$ running from
$u_{\pi(k)}$ to $I_{j_k}$, 
$k=1,2,\dots,p$, for some indices $j_{1},j_{2},\dots,j_p$ satisfying
$j_{1}<j_{2}<\dots<j_p$.

Then we have
$$
\sum _{\pi\in \frak S_{p}} ^{}(\sgn\pi)\cdot
\Cal P^{\text{nonint}}(\bold u_\pi\to I)
=
\Pf(Q),
\tag\AF
$$
with the matrix $Q=(Q_{i,j})_{1\le i,j\le p}$ given by 
$$Q_{i,j}=
\sum _{1\le u<v} ^{}
\big(\Cal P(u_i\to I_u)\cdot \Cal P(u_j\to I_v)-
\Cal P(u_j\to I_u)\cdot \Cal P(u_i\to I_v)\big),
\tag\AG
$$
where $\Cal P(A\to E)$ denotes the number of lattice paths from $A$ to
$E$.
\endproclaim

Application of this formula to our situation leads to the following
intermediate result.
In its statement, and throughout this paper, sums in which the lower
index is larger than the upper index have to be
interpreted according to the standard convention 
$$\sum _{r=m} ^{n-1}\text {\rm Expr}(k)=\cases \hphantom{-}
\sum _{r=m} ^{n-1} \text {\rm Expr}(k)&n>m\\
\hphantom{-}0&n=m\\
-\sum _{k=n} ^{m-1}\text {\rm Expr}(k)&n<m.\endcases
\tag\AH
$$

\proclaim{Proposition \TD}
Let $n,m,l$\/ be positive integers, and $k_1,k_2,\dots k_l$
positive integers with $k_1<k_2<\dots<k_l\le n/2$.
The number
of families 
$$(P_{-m+1},P_{-m+2},\dots,P_{m},P_{1^-},P_{2^-},\dots,P_{l^-},
P_{1^+},P_{2^+},\dots,P_{l^+})$$
of non-intersecting lattice paths,
where for $s\in\{-m+1,-m+2,\dots,m\}$, $P_s$ starts at 
$A_s=(s,-s+1)$,
while, 
for $t\in\{1,2,\dots,l\}$, $P_{t^-}$ starts at $B_{t^-}=(k_t,k_t+1)$ and 
$P_{t^+}$ starts at $B_{t^+}=(k_t+1,k_t)$, and all paths end somewhere
on the line $x+y=n+1$, is equal to
$$
(-1)^{\binom l2}\Pf(M),
\tag\AHa
$$
where $M$ is the skew-symmetric matrix with rows and columns indexed
by\linebreak 
$\{-m+1,-m+2,\dots,m,1^-,2^-,\dots,l^-,1^+,2^+,\dots,l^+\}$, and entries
given by
$$
M_{i,j}=\cases 
\sum_{r=i-j+1}^{j-i} \binom {2n}{n+r},&\text{if\/ }-m+1\le i<j\le m,\\
\sum_{r=i+1}^{-i} \binom {2n-2k_t}{n-k_t+r},&\text{if\/ $-m+1\le i\le
m$ and $j=t^-$},\\
\sum_{r=i}^{-i+1} \binom {2n-2k_t}{n-k_t+r},&\text{if\/ $-m+1\le i\le
m$ and $j=t^+$},\\
0,
&\text{if\/ $i=t^-$, $j={\hat t}^-$, and $1\le t<\hat t\le l$},\\
\binom {2n-2k_t-2k_{\hat t}}{n-k_t-k_{\hat t}}
+\binom {2n-2k_t-2k_{\hat t}}{n-k_t-k_{\hat t}+1},
&\text{if\/ $i=t^-$, $j={\hat t}^+$, and $1\le t,\hat t\le l$},\\
0,
&\text{if\/ $i=t^+$, $j={\hat t}^+$, and $1\le t<\hat t\le l$}.\\
\endcases
$$

\endproclaim

\demo{Proof}
We choose $p=2m+2l$, $u_i=A_{-m+i}$ for $i=1,2,\dots,2m$,
$u_{i}=B_{(i-2m)^-}$ for $i=2m+1,2m+2,\dots,2m+l$, 
$u_{i}=B_{(i-2m-l)^+}$ for $i=2m+l+1,2m+l+2,\dots,2m+2l$,
and $I_j=(j,n+1-j)$ for $j=-m+1,-m+2,\dots,n+m$
in Theorem~\TC. It is not difficult to convince oneself that, for
this choice of starting and ending points, all families of
nonintersecting lattice paths counted on the left-hand side of
(\AF) give rise to permutations $\pi$ which have the same sign,
and this sign is $(-1)^{\binom l2}$. Hence, the right-hand
side of (\AF) indeed counts the families of nonintersecting lattice
paths that we need to count, up to the overall sign $(-1)^{\binom l2}$. 
Then by Theorem~{\TC}, the number of these families is given 
by the Pfaffian $(-1)^{\binom l2}\Pf(Q)$, where the entries of the
skew-symmetric 
matrix $Q$ are given by (\AG). Consequently, for all pairs of indices
$(i,j)$ with $1\le i<j\le 2m+2l$, we have to compute the sum on the
right-hand side of (\AG). Since, in our situation, the definition of 
the starting points $u_i$ has three different cases, we have to
distinguish six cases. Since all of them are very similar, we treat
only one here in detail.

We choose the case of the $(i,j)$-entry of $Q$ with 
$1\le i<j\le 2m$. For this we have to compute the sum\footnote{ We
point out that the convention concerning row and column indices in
Theorem~\TC\ 
is different from the indexation which we use here.}
$$\multline 
\sum _{-m+1\le u<v\le n+m} ^{}
\Big(\Cal P\big((i,-i+1)\to (u,n+1-u)\big)\cdot 
\Cal P\big((j,-j+1)\to (v,n+1-v)\big)\\
-
\Cal P\big((j,-j+1)\to (u,n+1-u)\big)\cdot \Cal P\big((i,-i+1)\to
(v,n+1-v)\big)\Big). 
\endmultline$$
Using the simple fact that $\Cal P\big((a,b)\to (c,d)\big)=\binom
{c+d-a-b}{c-a}$, the above sum turns into
$$\align
&\sum _{-m+1\le u<v\le n+m} ^{}
\bigg(
\binom {n}{u-i}\binom n{v-j}
-
\binom n{u-j}\binom n{v-i}\bigg)\\
&\kern1cm
=
\sum _{v=0} ^{n+m}
\sum _{u=-m+1} ^{n+m}
\bigg(
\binom {n}{n-u+i}\binom n{u+v-j}
-
\binom n{n-u+j}\binom n{u+v-i}\bigg)\\
&\kern1cm
=
\sum _{v=0} ^{n+m}
\bigg(
\binom {2n}{n+v+i-j}
-
\binom {2n}{n+v+j-i}\bigg)\\
&\kern1cm
=
\sum _{v=0} ^{n+m}
\bigg(
\binom {2n}{n-v+j-i}
-
\binom {2n}{n-v+i-j}\bigg)\\
&\kern1cm
=\sum _{r=1} ^{j-i}
\binom {2n}{n+r}
+\sum _{r=i-j+1} ^{0}
\binom {2n}{n+r}
=\sum _{r=i-j+1} ^{j-i}
\binom {2n}{n+r},
\endalign$$
where we used the Chu--Vandermonde summation formula in the 
third line. This is exactly the corresponding expression in
the definition of $M_{i,j}$.\quad \quad \qed
\enddemo

Now we turn to the right-hand side of (\AC). We encode the lozenge
tilings of $H^-_{n,2m}(k_1,k_2,\dots,k_l)$ by families of
non-intersecting paths of lozenges (equivalently, lattice paths)
connecting vertical unit segments on its boundary.
See Figures~\FD\ and \FE\ for an example of this correspondence.
The power of 2 can then be absorbed by noticing that
$$
2^{n-2l}\,
\M^*\big(H^-_{n,2m}(k_1,k_2,\dots,k_l)\big)
$$
is equal to the weighted enumeration of
families $(P_1,P_2,\dots,P_m,P_{1^+},P_{2^+},\dots,P_{l^+})$ of 
non-intersecting lattice paths,
where, for $1\le s\le m$, the path $P_s$ runs from $A_s:=(s,-s+1)$
to $E_s:=(n+s,n-s+1)$, for $1\le t\le l$, the path $P_{t^+}$ 
runs from $B_{t^+}:=(k_t+1,k_t)$
to $E_s:=(n-k_t+1,n-k_t)$, the paths never cross the 
main diagonal $x=y$, 
and the weight of a family of paths is $2^{T}$, where $T$ is the number
of touching points of the paths in the family with the main diagonal.

Indeed, if there are $k$ such touching points, then there are
precisely $n-2l-k$ lozenges occupying positions weighted by $1/2$ (see
Figure~\FD) 
in the lozenge tiling of $H^-_{n,2m}(k_1,k_2,\dots,k_l)$ corresponding
to that family of lattice paths, and thus the weight of that tiling is  
$$
\left(\frac12\right)^{n-2l-k}=\left(\frac12\right)^{n-2l}2^k.
$$
Then the fraction on the right-hand side above cancels 
the factor of $2^{n-2l}$ on the right-hand side of (\AC), and the
above claim follows. 

By the classical Lindstr\"om--Gessel--Viennot formula for the enumeration
of non-intersecting lattice paths
(see \cite{\LindAA}\cite{\GeViAB}\cite{\StemAE}), the above weighted
count can be written 
in terms of a determinant. Below we recall this formula.

Let $G=(V,E)$ be a weighted directed acyclic graph with vertices $V$
and directed 
edges $E$, with weight function $w$ on its edges. The weight
$w(P)$ of a path $P$ in the graph 
is defined by $\prod _{e} ^{}w(e)$, where the
product is over all edges $e$ of the path.  
We denote the set of all paths in $G$ from $u$ to $v$ by $\P(u\to
v)$, and the set of all families $(P_1,P_2,\dots,P_p)$ of paths,
where $P_i$ runs from $u_i$ to $v_i$, $i=1,2,\dots,p$, by $\P(\bold
u\to\bold v)$, with $\bold u=(u_1,u_2,\dots,u_p)$ and 
$\bold v=(v_1,v_2,\dots,v_p)$. Denote by $\P^+(\bold u\to\bold v)$
the set of all families $(P_1,P_2,\dots,P_p)$ 
in $\P(\bold u\to\bold v)$ that are non-intersecting. 

The weight $w(\bold P)$ of a family $\bold
P=(P_1,P_2,\dots,P_p)$ of paths is defined as the product $\prod
_{i=1} ^{p}w(P_i)$ of all the weights of the paths in the family.
Finally, given a set $\Cal M$ with weight function $w$, we write 
$\GF(\Cal M;w)$ for the generating function $\sum
_{x\in\Cal M} ^{}w(x)$. 

\proclaim{Theorem \TE\ \smc(Lindstr\"om, Gessel and Viennot)}
With the above notation, we assume that the only permutation $\pi\in
S_p$ for which a family of non-intersecting lattice paths
$(P_1,P_2,\dots,P_p)$ exists such that the path $P_i$ connects
$u_i$ with $v_{\pi(i)}$ is the identity permutation. Then
$$
\GF(\P^+(\bold u\to \bold v);w)=
\det_{1\le i,j\le p}\big(\GF(\P(u_{j}\to v_i);w)\big).
\tag\AI$$
\endproclaim

Application of this formula to our situation leads to the following
intermediate result.

\topinsert
\vskip-60pt
\vbox{\noindent
\centertexdraw{
\drawdim truecm \setunitscale.7
\linewd.15
\move(0 0)
\RhombusA \RhombusA \RhombusB \RhombusA \RhombusA \RhombusB 
\RhombusA \RhombusB \RhombusB \RhombusA \RhombusA \RhombusA
\move(4.33013 -.5)
\RhombusB
\RhombusA \RhombusB \RhombusA \RhombusB \RhombusA 
\move(5.19615 0)
\RhombusB \RhombusA \RhombusB \RhombusA 
\move(4.33013 6.5)
\move(5.19615 7)
\move(6.06218 7.5)
\move(6.06218 .5)
\RhombusC 
\move(3.4641 0)
\RhombusA 
\move(6.9282 0)
\RhombusA 
\move(8.66025 0)
\RhombusA 
\move(12.12435 0)
\RhombusA 
\move(0 0)
\RhombusC \RhombusC 
\move(0 -1)
\RhombusC \RhombusC \RhombusC \RhombusC 
\move(0 -2)
\RhombusC \RhombusC \RhombusC \RhombusC \RhombusC 
\move(0 -3)
\RhombusC \RhombusC \RhombusC \RhombusC \RhombusC 
\move(3.4641 1)
\move(3.4641 0)
\RhombusC 
\move(2.59808 -.5)
\RhombusC \RhombusC \RhombusC 
\move(4.33013 -2.5)
\RhombusC \RhombusC 
\move(5.19615 -4)
\RhombusC \RhombusC 
\move(5.19615 -5)
\RhombusC \RhombusC 
\move(5.19615 1)
\move(6.9282 8)
\move(12.99038 -.5)
\RhombusB \RhombusB \RhombusB \RhombusB 
\move(9.52628  -.5)
\RhombusB \RhombusB 
\RhombusA \RhombusB  \RhombusB \RhombusB 
\move(7.794225 -.5)
\RhombusB \RhombusA \RhombusB  \RhombusB  
\RhombusA \RhombusB  \RhombusB 
\move(6.9282 1)
\move(6.9282 0)
\RhombusC 
\move(6.9282 -1)
\RhombusC \RhombusC 
\move(6.9282 -3)
\RhombusA \RhombusB \RhombusB \RhombusB \RhombusA
\move(6.9282 -5)
\RhombusB \RhombusB \RhombusA
\move(8.66025 -4)
\RhombusC
\move(8.66025 -5)
\RhombusC
\move(8.66025 0)
\RhombusC
\move(8.66025 1)
\move(8.66025 6)
\move(11.2583 3.5)
\move(8.66025 3)
\move(8.66025 3)
\move(10.3923 -1)
\RhombusB \RhombusA \RhombusB \RhombusB \RhombusB 
\move(11.2583 -.5)
\RhombusA \RhombusB \RhombusB \RhombusB \RhombusB 

\linewd.08
\move(0 -3.5)
\hdSchritt \hdSchritt \hdSchritt \hdSchritt \hdSchritt \odaSchritt
\hdSchritt \hdSchritt 
\odaSchritt \odaSchritt 
\hdSchritt \odaSchritt \odaSchritt \odaSchritt \odaSchritt \odaSchritt 
\move(0 -2.5)
\hdSchritt \hdSchritt \hdSchritt \hdSchritt \hdSchritt \odaSchritt
\hdSchritt \hdSchritt 
\odaSchritt \odaSchritt 
\hdSchritt \odaSchritt \odaSchritt \odaSchritt \odaSchritt \odaSchritt 
\move(0 -1.5)
\hdSchritt \hdSchritt \hdSchritt \hdSchritt \odaSchritt \hdSchritt 
\hdSchritt \odaSchritt 
\hdSchritt \odaSchritt 
\odaSchritt \hdSchritt 
\odaSchritt \odaSchritt \odaSchritt \odaSchritt 
\move(0 -.5)
\hdSchritt \hdSchritt \odaSchritt \hdSchritt \hdSchritt \hdSchritt 
\odaSchritt \odaSchritt 
\hdSchritt \hdSchritt 
\odaSchritt \odaSchritt \odaSchritt 
\hdSchritt \odaSchritt \odaSchritt 
\move(3.4641 -.5)
\hdSchritt \odaSchritt \odaSchritt 
\hdSchritt \hdSchritt \odaSchritt \hdSchritt \odaSchritt 

\ringerl(0 -3.5)
\ringerl(0 -2.5)
\ringerl(0 -1.5)
\ringerl(0 -.5)
\ringerl(13.8564 -3.5)
\ringerl(13.8564 -2.5)
\ringerl(13.8564 -1.5)
\ringerl(13.8564 -.5)
\ringerl(3.4641 -.5)
\ringerl(10.3923 -.5)

%
\rtext td:240 (3.35 -6.0){$\sideset {} \and {} \to 
    {\left.\vbox{\vskip2.4cm}\right\}}$}
\rtext td:0 (-1.3 -.2){$\sideset {2m} \and {}\to 
    {\left\{\vbox{\vskip2.4cm}\right.}$}
\htext (2.75 -7.0){$n$}
}
\vskip8pt
\centerline{\eightpoint A lozenge tiling of $H^-_{8,8}(2)$}
\vskip8pt
\centerline{\smc Figure \FD}
}
\vskip10pt
\endinsert

\topinsert
\vbox{\noindent
$$
\Gitter(14,10)(0,-4)
\Koordinatenachsen(14,10)(0,-4)
\Diagonale(-2,-2){12}
\Pfad(1,0),11211122 11222122\endPfad
\Pfad(2,-1),11112112 12212222\endPfad
\Pfad(3,-2),11111211 22122222\endPfad
\Pfad(4,-3),11111211 22122222\endPfad
\Pfad(3,2),12211212\endPfad
\DickPunkt(1,0)
\DickPunkt(2,-1)
\DickPunkt(3,-2)
\DickPunkt(4,-3)
\DickPunkt(3,2)
\DickPunkt(9,8)
\DickPunkt(10,7)
\DickPunkt(11,6)
\DickPunkt(12,5)
\DickPunkt(7,6)
\Kreis(0,0)
\Kreis(4,4)
\KREIS(4,4)
\Label\u{A_1}(1,0)
\Label\u{A_2}(2,-1)
\Label\u{A_3}(3,-2)
\Label\u{A_4}(4,-3)
\Label\o{B_{1^+}}(3,2)
\Label\ro{P_{1^+}}(5,3)
\Label\ro{P_4}(9,-4)
\Label\ro{P_3}(7,-2)
\Label\ro{P_2}(5,-1)
\Label\ro{P_1}(2,0)
\Label\u{x=y}(-2,0)
\Label\r{x}(14,0)
\Label\ru{y}(0,10)
\hskip6.5cm
$$
\centerline{\eightpoint The paths of lozenges of Figure {\FD} drawn as
non-intersecting lattice paths on $\Z^2$}
\vskip8pt
\centerline{\smc Figure \FE}
}
\vskip10pt
\endinsert

\proclaim{Proposition \TF}
Let $n,m,l$\/ be positive integers, and $k_1,k_2,\dots k_l$
positive integers with $k_1<k_2<\dots<k_l\le n/2$.
Then the generating function $\sum_{\bold P}w(\bold P)$ of all families 
$$\bold P=(P_1,P_2,\dots,P_m,P_{1^+},P_{2^+},\dots,P_{l^+})$$
of non-intersecting lattice paths,
where, for $1\le s\le m$, the path $P_s$ runs from $A_s=(s,-s+1)$
to $E_s=(n+s,n-s+1)$, for $1\le t\le l$, the path $P_{t^+}$ 
runs from $B_{t^+}=(k_t+1,k_t)$
to $E_t:=(n-k_t+1,n-k_t)$, all the paths never crossing over the 
main diagonal $x=y$, 
and the weight $w(\bold P)$ of a path family $\bold P$ 
being $2^{T(\bold P)}$ --- where $T(\bold P)$ is the number
of touching points of the paths in $\bold P$ with the main diagonal ---
is given by
$$
\det(N),
$$
where $N$ is the matrix with rows and columns indexed
by $\{1,2,\dots,m,1^+,2^+,\dots,l^+\}$, and entries
given by
$$
N_{i,j}=\cases 
\binom {2n}{n+j-i}+\binom {2n}{n-i-j+1},
&\text{if\/ }1\le i,j\le m,\\
\binom {2n-2k_t}{n-k_t-i+1}+\binom {2n-2k_t}{n-k_t-i},
&\text{if\/ $1\le i\le
m$ and $j=t^+$},\\
\binom {2n-2k_t}{n-k_t-j+1}+\binom {2n-2k_t}{n-k_t-j},
&\text{if\/ $i=t^+$ and $1\le j\le
m$ },\\
\binom {2n-2k_t-2k_{\hat t}}{n-k_t-k_{\hat t}}
+\binom {2n-2k_t-2k_{\hat t}}{n-k_t-k_{\hat t}-1},
&\text{if\/ $i=t^+$, $j={\hat t}^+$, and $1\le t,\hat t\le l$}.
\endcases
$$
\endproclaim

\demo{Proof}
We choose $p=m+l$, $u_i=A_{i}$ for $i=1,2,\dots,m$, and
$u_{i}=B_{(i-m)^+}$ for $i=m+1,m+2,\dots,m+l$
in Theorem~\TE. The underlying directed graph $G$ is
the graph whose vertices are the lattice points in $\Z^2$ lying
on or below the main diagonal $x=y$, and whose edges are the
horizontal edges $(x,y)\to(x+1,y)$ and vertical edges
$(x,y-1)\to(x,y)$ with $x\ge y$. The weight of horizontal edges is
$1$, as is the weight of vertical edges $(x,y-1)\to(x,y)$ for
$x>y$, while the weight of a vertical edge $(x,x-1)\to(x,x)$ is $2$.
Then it is not difficult to see that, for
this choice of starting and ending points, the technical
condition formulated at the beginning of the statement of
Theorem~\TE\ is satisfied. We can therefore apply the theorem.
In order to express the entries of the resulting determinant, we have
to compute the generating function 
$$
\GF\big((a,b)\to(c,d);w\big),
\tag\AJ
$$
where $a>b$ and $c>d$.
We claim that this generating function is given by
$$
\binom {c+d-a-b}{c-a}
+
\binom {c+d-a-b}{d-a}.
\tag\AK$$
Once this is shown,
the displayed expressions for the entries $N_{i,j}$ readily follow.

For the proof of our claim, we note that the generating function
in (\AJ) is a generating function for paths from $(a,b)$ to $(c,d)$
which never cross above the main diagonal $x=y$, and in which the
weight of a path $P$ is $2^{T(P)}$, where $T(P)$ is the number
of touching points of $P$ with the main diagonal. 
We may interpret this weight combinatorially as follows:
each of the above paths $P$ gives rise to $2^{T(P)}$ paths if
we reflect path portions between two successive touching points across the
main diagonal $x=y$, as well as the portion of the path from the
last touching point until the end point $(c,d)$. As a moment's thought
shows, in this way we obtain {\it all\/} paths from $(a,b)$ to
$(c,d)$ (meaning that we also obtain the ones which {\it do} cross over
above the main diagonal) as well as {\it all\/} paths from $(a,b)$ to
$(d,c)$. The total number of these paths is given by (\AK). This
completes the proof.\quad \quad \qed
\enddemo

\head 4. Equality of Pfaffian and determinant\endhead

We have seen in the previous section that,
in order to establish (\AC), we must show that the signed Pfaffian 
expression in (\AHa) equals the determinant in Proposition~\TF.
This is what we do in this section.

We should note that the matrix $M$ in Proposition~\TD\ has the form
$$
M=\pmatrix X&Y\\
-Y^t&Z\endpmatrix,
$$
where $X=(x_{j-i})_{-m+1\le i,j\le m}$ and
$Z=(z_{i,j})_{i,j\in\{1^-,\dots,l^-,1^+,\dots,l^+\}}$ 
are skew-symmetric, and $Y=(y_{i,j})_{-m+1\le i\le m,\,
j\in\{1^-,\dots,l^-,1^+,\dots,l^+\}}$ is a $2m\times 2l$ matrix.
Recalling the convention (\AH) of how to read sums, 
close inspection reveals that 
$
y_{i,t^-}
=-y_{-i,t^-}
$
and
$
y_{i,t^+}
=-y_{-i+2,t^+},
$
for all $i$ with $-m+1\le i\le m$ for which both sides of an equality are defined,
and $1\le t\le l$.
These properties are fundamental in the next lemma.

\proclaim{Lemma \TG}
For a positive integer $m$ and a non-negative integer $l$, let
$A$ be a matrix of the form
$$
A=\pmatrix X&Y\\
-Y^t&Z\endpmatrix,
$$
where $X=(x_{j-i})_{-m+1\le i,j\le m}$ and
$Z=(z_{i,j})_{i,j\in\{1^-,\dots,l^-,1^+,\dots,l^+\}}$ 
are skew-symmetric, 
and $Y=(y_{i,j})_{-m+1\le i\le m,\,
j\in\{1^-,\dots,l^-,1^+,\dots,l^+\}}$ is a $2m\times 2l$ matrix.
Suppose in addition that
$
y_{i,t^-}
=-y_{-i,t^-}
$
and
$
y_{i,t^+}
=-y_{-i+2,t^+},
$
for all $i$ with $-m+1\le i\le m$ for which both sides of an equality are defined, 
and $1\le t\le l$, and that
$z_{i,j}=0$ for all $i,j\in \{1^-,\dots,l^-\}$. Then
$$
\Pf(A)=(-1)^{\binom l2}\det(B),
$$
where 
$$
B=\pmatrix \bar X&\bar Y_1\\
\bar Y_2&\bar Z\endpmatrix,
$$
with 
$$\align 
\bar X&=(\bar x_{i,j})_{1\le i,j\le m},\\
\bar Y_1&=(y_{-i+1,j})_{1\le i\le m,\, j\in\{1^+,\dots,l^+\}}, \\
\bar Y_2&=(-y_{i,j})_{i\in\{1^-,\dots,l^-\},\,1\le j\le m}, \\
\bar Z&=(z_{i,j})_{i\in\{1^-,\dots,l^-\},\, j\in\{1^+,\dots,l^+\}}, 
\endalign$$
and the entries of $\bar X$ are defined by
$$
\bar x_{i,j}=x_{\vert j-i\vert+1}+x_{\vert j-i\vert+3}+\dots
+x_{i+j-1}.
$$
\endproclaim

\demo{Proof}
Starting from $A$, we construct a new matrix $\hat A$ as follows.
Replace the $i$-th row of $A$ by the sum
$$
\sum_{r=0}^{-i}(\text{row $i+2r$ of $A$}),
$$
for $i=0,-1,\dots,-m+1$.
In the resulting matrix, do the analogous operations with the
columns. Denote the matrix 
obtained this way by $\hat A$. Note that these operations do not
change the value of $\Pf(A)$, so $\Pf(A)=\Pf(\hat A)$. 

For $-m+1\le i,j\le 0$, the $(i,j)$-entry in $\hat A$ is
$$
\sum_{r=0}^{-i}\,\,\sum_{s=0}^{-j}
x_{j+2s-i-2r}
=
\sum_{t=i}^{-j} 
(\min\{-i,-j-t\}-\max\{0,-t\}+1)\,x_{j-i+2t}.
\tag\AL
$$
By assumption, the matrix $X$ is skew-symmetric, that is, we have
$x_r=-x_{-r}$ for $-m+1\le r\le m$. One readily verifies that
$$
\min\{-i,-j-t\}-\max\{0,-t\}+1
=\min\{-i,-j-(i-j-t)\}-\max\{0,-(i-j-t)\}+1,
$$
which is invariant under the replacement $t\to i-j-t$. Together with the
skew-symmetry,
this implies that the sum in (\AL)
vanishes, and, hence, the $(i,j)$-entry in $\hat A$ for $-m+1\le
i,j\le 0$ is zero.

Next we compute the $(i,j)$-entry of $\hat A$ for $-m+1\le i\le 0$
and $1\le j\le m$. This entry is not affected by the column
operations. Hence, it equals 
$$
\sum_{r=0}^{-i} x_{j-i-2r}
=\cases 
x_{j-i}+x_{j-i-2}+\dots+x_{j+i},
&\text{if $j+i\ge1$,}\\
x_{j-i}+x_{j-i-2}+\dots+x_{-j-i+2},
&\text{if $j+i\le0$,}
\endcases
\tag\AM
$$
where, in the second case, we used again that $x_r=-x_{-r}$ for all $r$.

Now we compute the $(i,j)$-entry of $\hat A$ for $-m+1\le i\le 0$
and $j=t^-$ with $1\le t\le l$. Again, this entry is not affected 
by the column operations, and, hence, it equals
$$
\sum_{r=0}^{-i} y_{i+2r,t^-}=0,
$$
where we used the property that $y_{r,t^-}=-y_{-r,t^-}$ for all $r$
and $t$. So also these entries vanish.

Finally, we compute the $(i,j)$-entry of $\hat A$ for $-m+1\le i\le 0$
and $j=t^+$ with $1\le t\le l$. Also this entry is not affected 
by the column operations, and, hence, it equals
$$
\sum_{r=0}^{-i} y_{i+2r,t^+}=
y_{i,t^+}+\sum_{r=1}^{-i} y_{i+2r,t^+}
=y_{i,t^+},
$$
where we used the property that $y_{r,t^+}=-y_{-r+2,t^+}$ for all $r$
and $t$. In other words, these entries do not change. The same is
true for all other entries above the diagonal which we did not yet
consider, as they are not affected by our row and column
operations. Since we did the same operation on the columns as we did
on the rows, the resulting matrix $\hat A$ is still skew-symmetric. 

Summarizing, the new matrix $\hat A$ has the block form
$$
\hat A=\pmatrix 0&\hat X&0&Y_{-,+}\\
-\hat X^t&X^+&Y_{+,-}&Y_{+,+}\\
0&-Y_{+,-}^t&0&Z_{-,+}\\
-Y_{-,+}^t&-Y_{+,+}^t&-Z_{-,+}^t&Z_{+,+}
\endpmatrix,
$$
where 
$$
\spreadlines{2\jot}
\align 
X^+&=(x_{j-i})_{1\le i,j\le m},\\
\hat X&=(\hat x_{i,j})_{-m+1\le i\le 0,\,1\le j\le m},
\endalign$$
with entry $\hat x_{i,j}$ given by (\AM), 
$$
\spreadlines{2\jot}
\align
Y_{-,+}&=(y_{i,j})_{-m+1\le i\le 0,\,j\in\{1^+,\dots,l^+\}},\\
Y_{+,-}&=(y_{i,j})_{1\le i\le m,\,j\in\{1^-,\dots,l^-\}},\\
Y_{+,+}&=(y_{i,j})_{1\le i\le m,\,j\in\{1^+,\dots,l^+\}},
\endalign
$$
and 
$$\align 
Z_{-,+}&=(z_{i,j})_{i\in\{1^-,\dots,l^-\},\,
j\in\{1^+,\dots,l^+\}},\\
Z_{+,+}&=(z_{i,j})_{i,j\in\{1^+,\dots,l^+\}}.
\endalign$$
By rearranging rows and columns in the same
way, this matrix may be brought into the form
$$
\pmatrix 
0&0&\hat X&Y_{-,+}\\
0&0&-Y_{+,-}^t&Z_{-,+}\\
-\hat X^t&Y_{+,-}&X^+&Y_{+,+}\\
-Y_{-,+}^t&-Z_{-,+}^t&-Y_{+,+}^t&Z_{+,+}
\endpmatrix.
$$
All these operations did not change the value of the Pfaffian,
except for a sign of $(-1)^{ml}$.
However, by the identity\footnote{ Up to sign, the identity follows by
using that the determinant is the square of the Pfaffian; the sign
follows by noticing that for $E=0$ this is Cayley's identity.} 
$$
\Pf\pmatrix 0&D\\-D^t&E\endpmatrix
=(-1)^{\binom d2}\det(D)
$$ 
for any $d\times d$ matrix $D$,
the Pfaffian of the last matrix is simply 
$$
\Pf(A)=\Pf(\hat A)=(-1)^{\binom {m+l}2+ml}
\det\pmatrix 
\hat X&Y_{-,+}\\
-Y_{+,-}^t&Z_{-,+}
\endpmatrix.
$$
Here, to be in line with the indexations of the blocks in this matrix,
the rows are indexed by
$i\in \{-m+1,-m+2,\dots,0,1^-,2^-,\dots,l^-\}$,
while the columns are indexed by 
$j\in \{1,2,\dots,m,1^+,2^+,\dots,l^+\}$.
We change the row index $i$, with $i=-m+1,-m+2,\dots,0$, 
to $-i+1$. 
This amounts to reversing the order of these rows.
The new row index will then range in $\{1,2,\dots,m,1^-,2^-,\dots,l^-\}$.
Hence, this leads to
$$\align
\Pf(A)=\Pf(\hat A)&=(-1)^{\binom m2+\binom {m+l}2+ml}
\det\pmatrix 
\bar X&\bar Y_1\\
-Y_{+,-}^t&Z_{-,+}
\endpmatrix\\
&=(-1)^{\binom l2}
\det\pmatrix 
\bar X&\bar Y_1\\
\bar Y_2&\bar Z
\endpmatrix,
\endalign$$
where $\bar X$, $\bar Y_1$, $\bar Y_2$, and $\bar Z$ are as in the
statement of the lemma. 
This is exactly $(-1)^{\binom l2}\det(B)$, and thus the proof is
complete.\quad \quad \qed 
\enddemo

\remark{Remark}
Gordon's Pfaffian reduction \cite{\GordAB, Lemma~1} is the special
case of Lemma~\TG\ where $l=0$. Indeed,
the row and column manipulations which we performed during the above proof
are the ones which Gordon used in his proof.
\endremark

\demo{Proof of Theorem \TA}
By Lemma~{\TG}, the signed Pfaffian in (\AHa)
is equal to
$$
\det\pmatrix 
x_{\vert j-i\vert+1}+x_{\vert j-i\vert+3}+\dots
+x_{i+j-1}&\sum_{r=-i+1}^{i} \binom {2n-2k_{\hat t}}{n-k_{\hat t}+r}\\
\sum_{r=-j+1}^{j} \binom {2n-2k_t}{n-k_t+r}&
\binom {2n-2k_t-2k_{\hat t}}{n-k_t-k_{\hat t}}
+\binom {2n-2k_t-2k_{\hat t}}{n-k_t-k_{\hat t}+1}\endpmatrix,
$$
where $1\le i,j\le m$ and $1\le t,\hat t\le l$, and 
$$
x_i=\sum_{r=-i+1}^{i} \binom {2n}{n+r}
$$
for all $i$.

The last series of operations consists in subtracting the $(m-1)$-st
row from the $m$-th, the $(m-2)$-nd row from the $(m-1)$-st, \dots,
the first from the second, and subsequently doing the analogous 
operations with the columns. One sees immediately that this converts
the above determinant into
$$
\det\pmatrix 
*
&
\binom {2n-2k_{\hat t}}{n-k_{\hat t}-i+1}+\binom {2n-2k_{\hat
t}}{n-k_{\hat t}+i}\\ 
\binom {2n-2k_t}{n-k_t-j+1}+\binom {2n-2k_t}{n-k_t+j}&
\binom {2n-2k_t-2k_{\hat t}}{n-k_t-k_{\hat t}}
+\binom {2n-2k_t-2k_{\hat t}}{n-k_t-k_{\hat t}+1}\endpmatrix,
$$
where the entries in the block marked by $*$ have still to be
computed. Comparison with the definition of the matrix $N$
in Proposition~\TF\ shows that the matrix of which the determinant is
taken above is precisely the same as $N$ in the top-right, bottom-left, and
bottom-right blocks. A somewhat more involved calculation 
shows that the top-left block (marked by $*$) of the above
matrix also agrees with the corresponding block of $N$.
This completes the proof of the equality of the signed Pfaffian (\AHa)
and the determinant in Proposition~\TF, and, hence, the proof of
Theorem~\TA.\quad \quad \qed
\enddemo

\remark{Remark}
The row and column operations which we performed during the above
proof are the ones which Stembridge uses to prove Theorem~7.1(a) in
\cite{\StemAE}. Again, the difference here is that our matrix has
$l$ more rows and $l$ more columns.
\endremark

\head 5. Odd holes\endhead

Since $k$ contiguous triangular holes of side two are equivalent, from
the point of view of counting the lozenge tilings of the region from
which they are removed, to a single triangular hole of side $2k$, one
sees that Theorem~{\TA} covers in fact the case when an arbitrary,
symmetric collection of even triangular holes is removed from along
the horizontal symmetry axis. What about odd holes? 

In order for the factorization identity to be meaningful, all three
kinds of tilings it involves need to exist. However, it is easy to
see that, if there is an odd
triangular hole with its vertical side not going through the center
of the underlying hexagon, 
there is no horizontally symmetric tiling. The only way for
odd holes to be present and horizontally symmetric tilings to exist is
if there are only two of them, symmetric about the center of the
hexagon and touching along their vertical edges (so as to form a
rhombus). In addition to them, we may have an arbitrary collection of
triangular holes of side two as in the previous sections. 

It turns out that the factorization identity holds in this more
general case as well. 

Let $n,m,l$\/ be positive integers, and let $x$ and
$k_1,k_2,\dotsc,k_l$ be positive integers so that 
$k_1<k_2<\dots<k_l\le n/2$. One readily sees that the region
$H_{n+x,2m}(k_1,k_2,\dotsc,\mathbreak k_l)$ has a horizontal lattice rhombus of
side length $x$ at its center precisely if $n$ is even. Assume that
this is so, and denote by $H_{n,2m}(k_1,k_2,\dotsc,k_l;x)$ the region
obtained from $H_{n+x,2m}(k_1,k_2,\dotsc,k_l)$ by removing from its
center this horizontal lattice rhombus of side $x$. 
(See Figure~\fea; in particular, we are still considering
the hexagon with side lengths $n+x,2m,n+x,n+x,2m,n+x$ 
with some holes inside.)

Then we have the following extension of Theorem~{\TA}.

\proclaim{Theorem \tea} For all positive integers $n,m,l$, with $n$
even, and all positive integers $x,k_1,k_2,\dotsc,k_l$ with 
$k_1<k_2<\dots<k_l\le n/2$, we have 
$$
\spreadlines{3\jot}
\multline
\M\left(H_{n,2m}(k_1,k_2,\dots,k_l;x)\right)\\
=
\M_{-}\left(H_{n,2m}(k_1,k_2,\dots,k_l;x)\right)
\M_{|}\left(H_{n,2m}(k_1,k_2,\dots,k_l;x)\right).
\endmultline
\tag\eea
$$

\endproclaim

\demo{Proof}
Note that for even $x$ this identity holds by Theorem~{\TA}. We deduce
it for odd values of $x$ by showing that, for fixed $n$ and $m$, and
fixed $k_1,\dotsc,k_l$, both sides of (\eea) are polynomial in $x$. 

Suppose, for the sake of the brevity of the write-up, that $l=0$,
i.e., there are no holes besides the two holes of side $x$ touching at
the center. We will see that polynomiality  in $x$ in the general case
follows by the same argument.

\topinsert
\centerline{\mypic{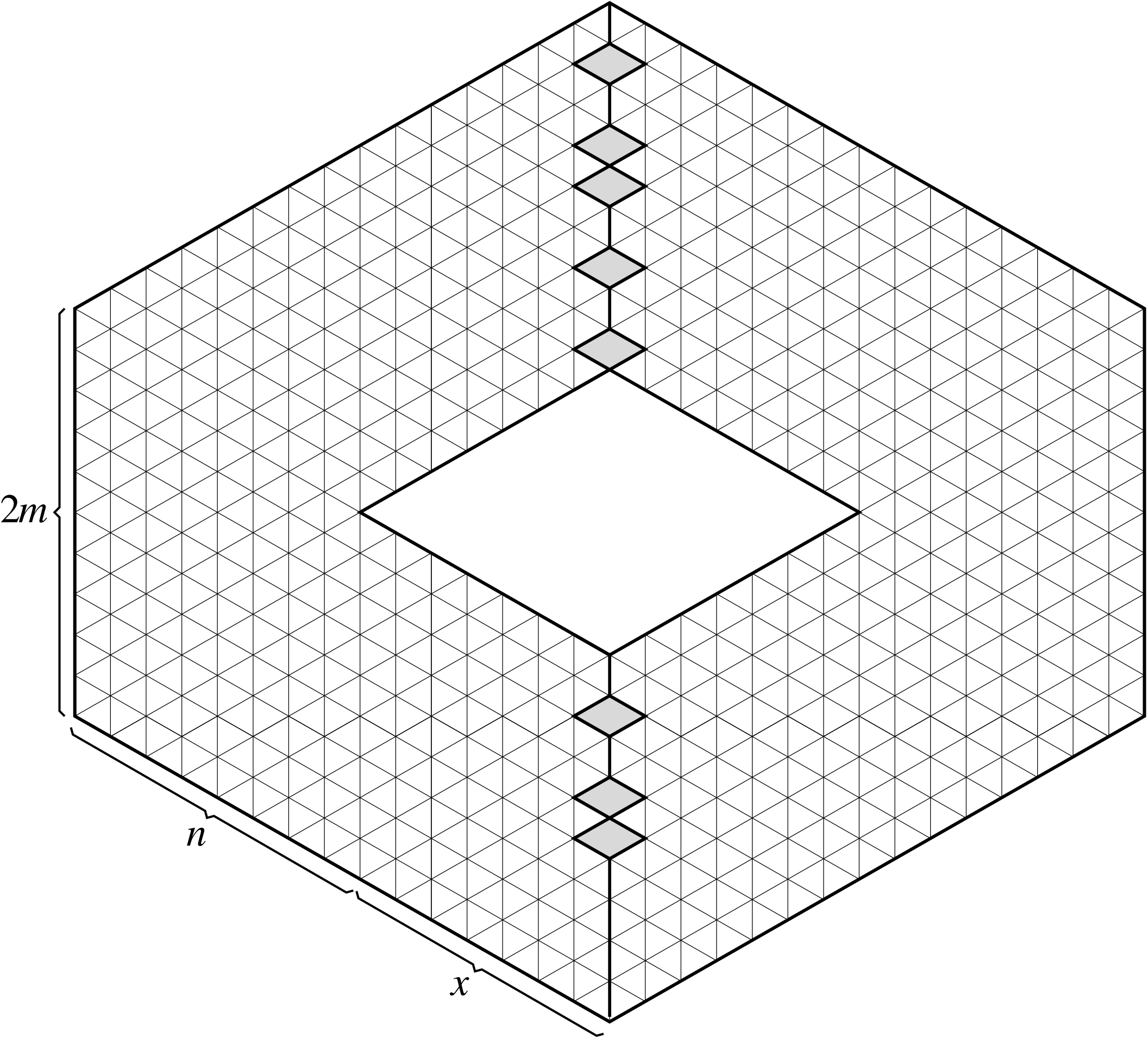}}
\centerline{{\smc Figure~\fea. {\rm The region $H_{15,10}(\emptyset;7)$.}}}
\endinsert

Let $m$ and $n$ be fixed, and allow $x$ to vary. In each tiling of
$H_{n,2m}(\emptyset;x)$, there are precisely $n$ horizontal lozenges
bisected by the vertical symmetry axis. Indeed, the tiling is encoded
by a family of $2m$ non-intersecting paths of lozenges running from the
left side of the hexagon to its right side, and each crosses the
vertical symmetry axis along a unit lattice segment; the remaining $n$
unit segments along that symmetry axis must therefore be occupied by
horizontal lozenges. 

Fix some set $S$ consisting of $n$ such horizontal lozenges, and
denote by $L_{n,2m}(S;x)$ (resp., $R_{n,2m}(S;x)$) the portion of
$H_{n,2m}(\emptyset;x)$ minus the union of the tiles in $S$ that is to
the left (resp., to the right) of the vertical symmetry axis. Then we
can write 
$$
\M\left(H_{n,2m}(\emptyset;x)\right)=
\sum_S \M\left(L_{n,2m}(S;x)\right) \M\left(R_{n,2m}(S;x)\right),
\tag\eeb
$$
where the sum ranges over all ${n+m\choose n}$ possible choices of the
set $S$.  

Since $L_{n,2m}(S;x)$ and $R_{n,2m}(S;x)$ are congruent, their number
of lozenge tilings is the same. Thus, by (\eeb), in order to prove
that $\M\left(H_{n,2m}(\emptyset;x)\right)$ is a polynomial in $x$, it
suffices to show that so is $\M\left(L_{n,2m}(S;x)\right)$, for any
$S$. 

This can be seen as follows. Encoding the lozenge tilings of
$L_{n,2m}(S;x)$ by paths of lozenges running from its left to
its right sides, and interpreting them as non-intersecting lattice
paths, we obtain that $\M\left(L_{n,2m}(S;x)\right)$ is equal to the
number of $2m$-tuples of non-intersecting lattice paths on $\Z^2$
having starting points $(2m-i,i-1)$, $i=1,\dotsc,2m$, and some fixed
$2m$-element subset of   
$$
\align
\!
\left\{(2m+n+x-j,j-1):j\in\left\{1,\dotsc,m+\frac{n}{2}\right\}
\cup\left\{m+\frac{n}{2}+x+1,\dotsc,2m+n+x\right\}\right\} 
\\
\tag\eec
\endalign
$$
as the set of its ending points. By the Lindstr\"om-Gessel-Viennot
theorem, this is given by a determinant of order $2m$, whose entries
are binomial coefficients. One readily sees that each of the resulting
binomial coefficients is a polynomial in $x$. It follows that so is
the determinant of the matrix, and the polynomiality of
$\M\left(H_{n,2m}(\emptyset;x)\right)$ as a function of $x$ follows. 

The very same argument proves the polynomiality in $x$ of
$\M_{|}\left(H_{n,2m}(\emptyset;x)\right)$. Indeed, the only
difference is that now the expression corresponding to (\eeb) only has
the first factor of the summand of (\eeb). 

The polynomiality in $x$ of $\M_{-}\left(H_{n,2m}(\emptyset;x)\right)$
follows by a similar argument, applied to the portion of
$H_{n,2m}(\emptyset;x)$ that is above a zig-zag lattice line analogous
to the one in Figure~{\fbb}. 

Consider now the case when $l$ is not necessarily zero. The effect of
the presence of the $l$ holes of side two in the left half of the
hexagon is to introduce $2l$ more lattice paths (except when looking
at the horizontally symmetric tilings, when it is only $l$ more
lattice paths), with starting points having some fixed coordinates
(not involving $x$), and ending points still chosen from among the
elements of the set (\eec). It is then clear that the above arguments
work equally well in this general case, showing that all three
quantities in (\eea) are polynomials in $x$. Since the two sides of
(\eea) agree for all even values  of $x$ (cf\. Theorem~\TA), it follows
that they agree for all $x$. This completes the proof.\quad \quad \qed
\enddemo

\remark{Remark}
The case $l=0$, $x=1$ was first proved, by a different method, by
Kasraoui and the second author in \cite{\KK}. An alternative way
to view this, is that a combination of this case of Theorem~\tea\
with the main result in \cite{\CiKrAA} yields a new proof of the
main result in \cite{\KK}.
\endremark

\head 6. Concluding remarks\endhead

The factorization result proved in this paper is reminiscent of the
symmetries considered by Kuperberg in \cite{\KupExploration,
Sec.~IV.C}, especially given that, in the terminology
of \cite{\KupExploration}, one of our symmetries is color-preserving
and the other color-reversing. However, it turns out that the
factorization we prove in this paper does not  hold even under small
changes in the structure of the holes we considered. The fact that our
factorization result depends so strongly on the geometry of the holes
indicates that it is of a kind different from the ones considered
in \cite{\KupExploration}. 

Given the equivalence of lozenge tilings of a hexagon with
plane partitions contained in a given box (cf\. \cite{\DT}),
one way to view the special case $l=0$ of Theorem~{\TA} is that it
gives a new proof for the enumeration of symmetric plane
partitions contained in a given box 
(first proved by Andrews \cite{\AndrAK}).
Indeed, by our factorization result, the latter is the ratio
of the total number of plane partitions that fit in the box
and the number of those that are transpose-complementary.
The formula for the total number (due to MacMahon \cite{\MacMAA}) can easily
be proved inductively using Kuo's graphical condensation (cf\.
\cite{\Kuo}), while the transpose-complementary case 
(first proved by Proctor \cite{\Proc}) can be directly deduced from
MacMahon's result by applying the matchings factorization
theorem (cf\. \cite{\CiucAB}). From
the point of 
view of \cite{\CiKrfiveandahalf}, this adds a seventh symmetry class
that can be proved in a combinatorial way. 

It would be interesting to find some more direct relations between the
symmetry classes of plane partitions, and also to find an algebraic
generalization of the factorization presented in this paper in the
style of \cite{\CiKrAF}, which generalized the special case $l=0$ of
the main theorem here in terms of Schur functions.

\enddocument